\documentclass[aos,MSNbibl,dvips]{arximspdf}
\usepackage{subenv}
\usepackage{graphicx}

\doi{10.1214/12-AOS1033} 
\volume{40}
\issue{5}
\pubyear{2012}
\firstpage{2483}
\lastpage{2510}

\makeatletter

\newcommand{\rright}{\right}
\newcommand{\lleft}{\left}

\newtheorem{theorem}{Theorem}
\newtheorem{corollary}{Corollary}
\newtheorem{lemma}{Lemma}

\newproclaim{example}{Example}
\newproclaim{Remarks}{Remarks}

\newcommand{\llangle}{\langle\!\langle}
\newcommand{\rrangle}{\rangle\!\rangle}

\newcommand{\llanglel}{\bigl\langle\bigl\langle}
\newcommand{\rrangler}{\bigr\rangle\bigr\rangle}

\newcommand{\nnorm}{|\!|\!|}
\newcommand{\nnnorml}{\bigl|\!\bigl|\!\bigl|}
\newcommand{\nnnormr}{\bigr|\!\bigr|\!\bigr|}

\newcommand{\wnoise}[1]{w_{#1}}
\newcommand{\yobs}[1]{y_{#1}}

\newcommand{\Hil}{\mathcal{H}}
\newcommand{\mprob}{\mathbb{P}}

\newcommand{\ct}{\mathcal{T}}
\newcommand{\reals}{\mathbb{R}}
\newcommand{\real}{\mathbb{R}}
\newcommand{\widebar}{\overline}

\newcommand{\Exs}{\mathbb{E}}
\newcommand{\diag}{\operatorname{diag}}
\newcommand{\trace}{\operatorname{tr}}

\newcommand{\xobs}[1]{x_{#1}}
\newcommand{\Kbb}{\mathbb{K}}
\newcommand{\dist}{\mathrm{d}}
\newcommand{\Null}{\operatorname{Ker}}

\newcommand{\Span}{\operatorname{span}}
\newcommand{\rp}{R_{m}}

\newcommand{\MiniMax}{\mathcal{M}_{m, n}}
\newcommand{\Ball}{\mathbb{B}}

\newcommand{\epscrit}{\epsilon_{m, n}}
\newcommand{\SimpleF}{\mathcal{F}}
\newcommand{\AnnoyCon}{\kappa_{r, \rho}}
\newcommand{\myspan}{\operatorname{span}}

\newcommand{\Range}{\operatorname{Ra}}
\newcommand{\Ker}{\operatorname{Ker}}
\newcommand{\KerFun}{\mathbb{K}}

\newcommand{\colsp}{\operatorname{col}}

\makeatother

\begin{document}
\begin{frontmatter}

\title{Sampled forms of functional PCA in reproducing kernel
Hilbert spaces}
\runtitle{Sampled form of functional PCA in RKHS}

\begin{aug}
\author[A]{\fnms{Arash A.} \snm{Amini}\corref{}\ead[label=e1]{amini@eecs.berkeley.edu}}
\and
\author[B]{\fnms{Martin J.} \snm{Wainwright}\thanksref{T1}\ead[label=e2]{wainwrig@stat.berkeley.edu}}
\runauthor{A. A. Amini and M. J. Wainwright}
\affiliation{University of California, Berkeley}
\address[A]{Department of Electrical Engineering \\
\quad and Computer Science \\
University of California, Berkeley \\
Berkeley, California 94720\\
USA\\
\printead{e1}}
\address[B]{Department of Statistics \\
University of California, Berkeley \\
Berkeley, California 94720\\
USA\\
\printead{e2}} 
\end{aug}

\thankstext{T1}{Supported in part by NSF Grants
CCF-0545862 and DMS-09-07632.}

\received{\smonth{9} \syear{2011}}
\revised{\smonth{4} \syear{2012}}

%
\begin{abstract}
We consider the sampling problem for functional PCA (fPCA), where
the simplest example is the case of taking time samples of the
underlying functional components. More generally, we model the
sampling operation as a continuous linear map from $\Hil$ to
$\reals^m$, where the functional components to lie in some
Hilbert subspace $\Hil$ of $L^2$, such as a reproducing kernel
Hilbert space of smooth functions. This model includes time and
frequency sampling as special cases. In contrast to classical
approach in fPCA in which access to entire functions is assumed,
having a limited number $m$ of functional samples places
limitations on the performance of statistical procedures. We study
these effects by analyzing the rate of convergence of an
$M$-estimator for the subspace spanned by the leading components in
a multi-spiked covariance model. The estimator takes the form of
regularized PCA, and hence is computationally attractive. We
analyze the behavior of this estimator within a nonasymptotic
framework, and provide bounds that hold with high probability as a
function of the number of statistical samples
$n$ and the number of functional samples $m$. We also
derive lower bounds showing that the rates obtained are minimax
optimal.
\end{abstract}

%

\end{frontmatter}

\section{Introduction}
The statistical analysis of functional data, commonly referred to as
functional data analysis (FDA), is an established area of statistics
with a great number of practical applications; see the
books~\cite{Ramsay2005,Ramsay2002} and references therein for various
examples. When the data is available as finely sampled curves, say in
time, it is common to treat it as a collection of continuous-time
curves or functions, each being observed in totality. These datasets
are then termed ``functional,'' and various statistical procedures
applicable in finite dimensions can be extended to this functional
setting. Among such procedures is principal component analysis (PCA),
which is the focus of present work.

If one thinks of continuity as a mathematical abstraction of reality,
then treating functional data as continuous curves is arguably a valid
modeling device. However, in practice, one is faced with finite
computational resources and is forced to implement a
(finite-dimensional) approximation of true functional procedures by
some sort of truncation procedure, for instance, in the frequency
domain. It is then important to understand the effects of this
truncation on the statistical performance of the procedure. In other
situations, such as in longitudinal data analysis~\cite{Diggle2002},
a~continuous curve model is justified as a hidden underlying generating
process to which one has access only through sparsely sampled
measurements in time, possibly corrupted by noise. Studying how the
time-sampling affects the estimation of the underlying functions in the
presence of noise shares various common elements with the
frequency-domain problem described above.

The aim of this paper is to study effects of ``sampling''---in a
fairly general sense---on functional principal component analysis in
smooth function spaces. In order to do so, we adopt a
functional-theoretic approach by treating the sampling procedure as a
(continuous) linear operator. This set-up provides us with a notion of
sampling general enough to treat both the frequency-truncation and
time-sampling within a unified framework. We take as our smooth
function space a Hilbert subspace $\Hil$ of $L^2[0,1]$ and denote the
sampling operator by $\Phi\dvtx  \Hil
\to\reals^m$. We assume that there are functions $\xobs{i}(t)$, $t
\in[0,1]$, in $\Hil$ for $i=1,\ldots,n$, generated i.i.d. from a
probabilistic model (to be discussed). We then observe the collection
$\{ \Phi x_i \}_{i=1}^n\subset\reals^m$ in noise. We
refer to the index $n$ as the number of \textit{statistical
samples}, and to the index $m$ as the number of \textit{functional
samples}.

We analyze a natural $M$-estimator which takes the form of a
regularized PCA in~$\reals^m$, and provide nonasymptotic bounds
on the estimation error in terms of $n$ and~$m$. The
eigen-decay of two operators govern the rates, the product of the
sampling operator $\Phi$ and its adjoint, and the product of the map
embedding $\Hil$ in $L^2$ and its adjoint. Our focus will be on the
setting where $\Hil$ is a reproducing kernel Hilbert space (RKHS), in
which case the two eigen-decays are intimately related through the
kernel function $(s,t) \mapsto
\Kbb(s,t)$. In such cases, the two components of the rate interact and
give rise to optimal values for the number of functional samples
($m$) in terms of the number of statistical samples ($n$) or
vice versa. This has practical appeal in cases where obtaining either
type of samples is costly.

Our model for the functions $\{\xobs{i}\}$ is an extension to function
spaces of the \textit{spiked covariance model} introduced by Johnstone
and his collaborators~\cite{John01,Johnstone2004}, and studied by
various authors (e.g.,~\cite{Johnstone2004,Paul2008,AmiWai2009}). We
consider such models with $r$ components, each lying within the
Hilbert ball $\Ball_\Hil(\rho)$ of radius $\rho$, with the goal of
recovering the $r$-dimensional subspace spanned by the spiked
components in this functional model. We analyze our $M$-estimators
within a high-dimensional framework that allows both the number of
statistical samples $n$ and the number of functional samples
$m$ to diverge together. Our main theoretical contributions are
to derive nonasymptotic bounds on the estimation error as a function
of the pair $(m, n)$, which are shown to be sharp
(minimax-optimal). Although our rates also explicitly track the
number of components $r$ and the smoothness parameter $\rho$, we
do not make any effort to obtain optimal dependence on these
parameters.

The general asymptotic properties of PCA in function spaces have been
investigated by various authors
(e.g.,~\cite{Dauxois1982,Bosq2000,Hall2006a}). Accounting for
smoothness of functions by introducing various roughness/smoothness
penalties is a standard approach, used in the
papers~\cite{Rice1991,Pezzulli1993,Silverman1996,Boente2000}, among
others. The problem of principal component analysis for sampled
functions, with a similar functional-theoretic perspective, is
discussed by Besse and Ramsey~\cite{Besse1986} for the noiseless case.
A more recent line of work is devoted to the case of functional PCA
with noisy sampled
functions~\cite{Cardot2000,Yao2005,Hall2006}. Cardot~\cite{Cardot2000}
considers estimation via spline-based approximation, and derives MISE
rates in terms of various parameters of the model. Hall et
al.~\cite{Hall2006} study estimation via local linear smoothing, and
establish minimax-optimality in certain settings that involve a fixed
number of functional samples. Both papers~\cite{Cardot2000,Hall2006}
demonstrate trade-offs between the numbers of statistical and
functional samples; we refer the reader to Hall et al.~\cite{Hall2006}
for an illuminating discussion of connections between FDA and LDA
approaches (i.e., having full versus sampled functions), which inspired
much of the present work. We note that the regularization present in
our $M$-estimator is closely related to classical roughness
penalties~\cite{Rice1991,Silverman1996} in the special case of spline
kernels, although the discussion there applies to fully-observed
functions, as opposed to the sampled models considered here.

After initial posting of this work, we became aware of more
recent work on sampled functional PCA. Working within the framework
of Hall et al.~\cite{Hall2006}, the analysis of Li and
Hsing~\cite{LiHsi10} allows for more flexible sample sizes per curve;
they derive optimal uniform (i.e., $L^\infty$) rates of convergence
for local linear smoothing estimators of covariance function and the
resulting eigenfunctions. Another line of work~\cite{HuaETA08,QiZha11}
has analyzed sampled forms of Silverman's
criterion~\cite{Silverman1996}, with some variations. Huang et
al.~\cite{HuaETA08} derive a criterion based on rank-one approximation
coupled with scale invariance considerations, combined with an extra
weighting of the covariance matrix. Xi and Zhao~\cite{QiZha11} also
show the consistency of their estimator for both regular and irregular
sampling. The regular (time) sampling setup in both papers have an
overlap with our work; the eigenfunctions are assumed to lie in a
second order Sobolev space, corresponding to a special case of a RKHS.
However, even in this particular case, our estimator is different, and
it is an interesting question whether a version of the results presented
here can be used to show the minimax optimality of these
Silverman-type criteria. There has also been recent work with emphasis
on sampled functional covariance estimation, including the work of Cai
and Yuan~\cite{CaiYua10}, who analyze an estimator which can be
described as regularized least-squares with penalty being the norm of
tensor product of RKHS with itself. They provide rates of convergence
for the covariance function, from which certain rates (argued to be
optimal within logarithmic factors) for eigenfunctions follow.

As mentioned above, our sampled model resembles very much that of
spiked covariance model for high-dimensional principal component
analysis. A line of work on this model has treated various types of
sparsity conditions on the
eigenfunctions~\cite{Johnstone2004,Paul2008,AmiWai2009}; in contrast,
here the smoothness condition on functional components translates into
an ellipsoid condition on the vector principal components. Perhaps an
even more significant difference is that in this paper, the effective
scaling of noise in $\reals^m$ is substantially smaller in some
cases (e.g., the case of time sampling).
This difference could explain why the difficulty of
``high-dimensional'' setting is not observed in such cases as one lets
$m,n\to\infty$. On the other hand, a difficulty
particular to our sampled model is the lack of orthonormality between
components after sampling. It not only leads to identifiability issues, but
also makes recovering individual components difficult.

In order to derive nonasymptotic bounds on our $M$-estimator, we
exploit various techniques from empirical process theory
(e.g.,~\cite{Geer2009}), as well as the concentration of measure
(e.g.,~\cite{Ledoux2001}). We also exploit recent
work~\cite{Mendelson2002} on the localized Rademacher complexities of
unit balls in a reproducing kernel Hilbert space, as well as
techniques from nonasymptotic random matrix theory, as discussed in
Davidson and Szarek~\cite{Davidson2001}, in order to control various
norms of random matrices. These techniques allow us to obtain
finite-sample bounds that hold with high probability, and are
specified explicitly in terms of the pair $(m, n)$, and the
underlying smoothness of the Hilbert space.

The remainder of this paper is organized as follows.
Section~\ref{SecBackground} is devoted to background material on
reproducing kernel Hilbert spaces, adjoints of operators, as well as
the class of sampled functional models that we study in this paper.
In Section~\ref{SecMest}, we describe $M$-estimators for sampled
functional PCA, and discuss various implementation details.
Section~\ref{SecMain} is devoted to the statements of our main
results, and discussion of their consequences for particular sampling
models. In subsequent sections, we provide the proofs of our results,
with some more technical aspects deferred to the supplementary material~\cite{supp}.
Section~\ref{SecProofSubspace} is devoted to bounds on the
subspace-based error. We conclude with a discussion in
Section~\ref{SecDiscuss}. In the supplementary material~\cite{supp}, Section 7 
is devoted to
proofs of bounds on error in the function space, whereas Section 8 %
provides proofs of matching lower bounds on the minimax error, showing
that our
analysis is sharp.

\textit{Notation}. We will use \mbox{$\nnorm\cdot\nnorm_{\mathrm{HS}}$} to denote the
Hilbert--Schmidt norm of an operator or a matrix. The corresponding
inner product is denoted as $\llangle{\cdot},{\cdot}\rrangle$. If
$T$ is an
operator on a Hilbert space $\Hil$ with an orthonormal basis
$\{e_j\}$, then $\nnorm T \nnorm_{\mathrm{HS}}^2 = \sum_{j} \|T e_j\|_{\Hil
}^2$. For a
matrix $A = (a_{ij})$, we have $\nnorm A \nnorm_{\mathrm{HS}}^2 = \sum_{i,j}
|a_{ij}|^2$. For a linear operator $\Phi$, the adjoint is
denoted as $\Phi^*$, the range as $\Range(\Phi)$ and the kernel as
$\Ker(\Phi)$.


\section{Background and problem set-up}
\label{SecBackground}

In this section, we begin by introducing background on reproducing
kernel Hilbert spaces, as well as linear operators and their adjoints.
We then introduce the functional and observation model that we study
in this paper, and conclude with discussion of some
approximation-theoretic issues that play an important role in parts of
our analysis.

\subsection{Reproducing kernel Hilbert spaces}
\label{SecRKHS}

We begin with a quick overview of some standard properties of
reproducing kernel Hilbert spaces; we refer the reader to the
books~\cite{Wahba,Gu02} and references therein for more details.
A~reproducing kernel Hilbert space (or RKHS for short) is a Hilbert
space $\Hil$ of functions $f\dvtx
T\to\reals$ that is equipped with a symmetric positive semidefinite
function $\Kbb\dvtx
T\times T\to\reals$, known as the kernel function.
We assume the kernel to be continuous, and the set
$T\subset\reals^d$ to be compact. For concreteness, we think
of $T= [0,1]$ throughout this paper, but any compact set of
$\real^d$ suffices. For each $t \in T$, the function
$R_t:=\Kbb(\cdot, t)$ belongs to the Hilbert space $\Hil$ and it
acts as the \textit{representer of evaluation}, meaning that
$\langle f,R_t\rangle_{\Hil} = f(t)$ for all $f \in\Hil$.

The kernel $\Kbb$ defines an integral operator $\ct_\Kbb$ on
$L^2(T)$, mapping the function $f$ to the function $g(s) =
\int_T\Kbb(s,t) f(t) \,dt$. By the spectral theorem in Hilbert
spaces, this operator can be associated with a sequence of
eigenfunctions $\psi_k, k=1,2,\ldots\,$, in $\Hil$, orthogonal in
$\Hil$ and orthonormal in $L^2(T)$, and a sequence of
nonnegative eigenvalues $\mu_1 \geq\mu_2 \geq\cdots\,$. Most useful
for this paper is the fact that any function $f \in\Hil$ has an
expansion in terms of these eigenfunctions and eigenvalues, namely
%
\begin{equation}
\label{eqfalphaisomorphism} f = \sum_{k=1}^\infty
\sqrt{\mu_k} \alpha_k \psi_k
\end{equation}
for some $(\alpha_k) \in\ell^2$. In terms of this expansion, we have
the representations $\| f \|_\Hil^2 = \sum_{k=1}^\infty
\alpha_k^2$ and $\| f \|_{L^2}^2 = \sum_{k=1}^\infty\mu_k
\alpha_k^2$. Many\vspace*{1pt} of our results involve the decay rate of these
eigenvalues: in particular, for some parameter $\alpha> 1/2$, we say
that the kernel operator has eigenvalues with
\textit{polynomial-$\alpha$ decay} if there is a constant $c > 0$ such
that
%
\begin{equation}
\label{EqnPolySmooth} \mu_k \leq\frac{c}{k^{2 \alpha}} \qquad\mbox{for
all $k = 1, 2, \ldots\,$.}
\end{equation}
Let us consider an example to illustrate.
%
\begin{example}[(Sobolev class with smoothness $\alpha= 1$)]
In the case
$T= [0,1]$ and $\alpha= 1$, we can consider the
kernel function $\Kbb(s,t) = \min\{s, t \}$. As discussed in
Appendix A 
of the supplementary material~\cite{supp}, this kernel generates the
class of functions
\[
\Hil:= \bigl\{ f \in L^2\bigl([0,1]\bigr) \mid f(0) =0, f \mbox{
absolutely continuous and } f' \in L^2\bigl([0,1]\bigr)
\bigr\}.
\]
The class $\Hil$ is an RKHS with inner product $\langle f,g\rangle_\Hil=
\int_0^1 f'(t) g'(t) \,dt$, and the ball $\Ball_\Hil(\rho)$ corresponds
to a Sobolev space with smoothness $\alpha= 1$. The
eigen-decomposition of the kernel integral operator is
%
\begin{equation}\quad
\mu_k = \biggl[ \frac{(2k-1) \pi}{2} \biggr]^{-2},\qquad
\psi_k(t) = \sqrt{2} \sin\bigl( \mu_k^{-1/2} t
\bigr),\qquad k=1,2,\ldots.
\end{equation}
Consequently, this class has polynomial decay with parameter $\alpha= 1$.
\end{example}
We note that there are natural generalizations of this
example to $\alpha= 2, 3, \ldots\,$, corresponding to the
Sobolev classes of $\alpha$-times differentiable functions; for
example, see
the books~\cite{Wahba,Gu02,BerTho04}.

In this paper, the operation of generalized sampling is defined in
terms of a bounded linear operator $\Phi\dvtx  \Hil\rightarrow
\real^m$ on the Hilbert space. Its adjoint is a mapping
$\Phi^*\dvtx \real^m\rightarrow\Hil$, defined by the relation
$\langle\Phi f,a\rangle_{\real^m} = \langle f,\Phi^* a\rangle_\Hil$ for
all $f \in\Hil$ and $a \in\real^m$. In order to compute a
representation of the adjoint, we note that by the Riesz
representation theorem, the $j$th coordinate of this
mapping---namely, $f \mapsto[\Phi f]_j$---can be represented as an
inner product $\langle\phi_j,f\rangle_\Hil$, for some element \mbox{$\phi_j
\in
\Hil$}, and we can write
%
\begin{equation}
\label{EqnLinopRep} \Phi f = \lleft[\matrix{ \langle\phi_1,f
\rangle_\Hil& \langle\phi_2,f\rangle_\Hil&
\cdots& \langle\phi_m,f\rangle_\Hil} \rright]^T.
\end{equation}
Consequently,\vspace*{1pt} we have $\langle\Phi f,a\rangle_{\real^m} =
\sum_{j=1}^ma_j \langle\phi_j,f\rangle_\Hil =
\langle\sum_{j=1}^ma_j \phi_j,f\rangle_\Hil$, so that for any
$a \in\real^m$, 
%
\begin{equation}
\label{EqnAdjoint} \Phi^*a = \sum_{j=1}^ma_j
\phi_j.
\end{equation}
This adjoint operator plays an important role in our analysis.


\subsection{Functional model and observations}

Let $s_1 \geq s_2 \geq s_3 \geq\cdots\geq s_r
> 0$ be a
fixed sequence of positive numbers, and let $\{f^*_j\}_{j
=1}^r$ be a fixed sequence of functions orthonormal in
$L^2[0,1]$. Consider a collection of $n$ i.i.d. random
functions $\{x_1,\ldots, x_n\}$, generated according to the
model
%
\begin{equation}
\label{EqnXmodel} \xobs{i}(t) = \sum_{j=1}^rs_j
\beta_{ij} f^*_j(t) \qquad\mbox{for $i = 1,\ldots,n$},
\end{equation}
where $\{\beta_{ij}\}$ are i.i.d. $N(0,1)$ across all pairs $(i,j)$.
This model corresponds to a finite-rank instantiation of functional
PCA,\vadjust{\goodbreak} in which the goal is to estimate the span of the unknown
eigenfunctions $\{f^*_j\}_{j=1}^r$. Typically, these
eigenfunctions\vspace*{1pt} are assumed to satisfy certain smoothness
conditions; in this paper, we model such conditions by assuming that
the eigenfunctions\vspace*{1pt} belong to a reproducing kernel Hilbert space $\Hil$
embedded within $L^2[0,1]$; more specifically, they lie in some ball in
$\Hil$,
%
\begin{equation}
\label{eqhilballassump} \| f^*_j \|_\Hil\le\rho,\qquad j=1,
\ldots,r.
\end{equation}

For statistical problems involving estimation of functions, the random
functions might only be observed at certain times $(t_1,\ldots,
t_m)$, such as in longitudinal data analysis, or we might collect
only projections of each $\xobs{i}$ in certain directions, such as in
tomographic reconstruction. More concretely, in a \textit{time-sampling
model}, we observe $m$-dimensional vectors of the form
%
\begin{equation}
\label{EqnTimeSamp} y_i = \lleft[\matrix{
\xobs{i}(t_1) & \xobs{i}(t_2) & \cdots& \xobs
{i}(t_m) } \rright]^T + \sigma_0
w_i\qquad\mbox{for $i = 1, 2,\ldots, n$},
\end{equation}
where $\{t_1, t_2,\ldots, t_m\}$ is a fixed collection of design
points, and $w_i \in\real^m$ is a noise vector. Another
observation model is the \textit{basis truncation model} in which we
observe the projections of $f$ onto the first $m$ basis functions
$\{\psi_j\}_{j=1}^m$ of the kernel operator---namely,
%
\begin{eqnarray}
\label{EqnBasisTrun} y_i = \lleft[\matrix{ \langle
\psi_1,\xobs{i}\rangle_{L^{2}} & \langle\psi_2,
\xobs{i}\rangle_{L^{2}} & \cdots& \langle\psi_m,\xobs{i}
\rangle_{L^{2}} } \rright]^T + \sigma_0
w_i\nonumber\\[-8pt]\\[-8pt]
&&\eqntext{\mbox{for $i = 1, 2,\ldots, n$},}
\end{eqnarray}
where $\langle\cdot,\cdot\rangle_{L^{2}}$ represents the inner
product in
$L^2[0,1]$.

In order to model these and other scenarios in a unified manner, we
introduce a linear operator $\Phi_m$ that maps any function $x$
in the Hilbert space to a vector $\Phi_m(x)$ of $m$ samples,
and then consider the linear observation model
%
\begin{equation}
\label{EqnLinObs} y_i = \Phi_m(\xobs{i}) +
\sigma_mw_i\qquad\mbox{for $i = 1, 2,\ldots, n$.}
\end{equation}
This model (\ref{EqnLinObs}) can be viewed as a functional analog of
the spiked covariance models introduced by
Johnstone~\cite{John01,Johnstone2004} as an analytically-convenient
model for studying high-dimensional effects in classical PCA.

Both the time-sampling (\ref{EqnTimeSamp}) and frequency
truncation (\ref{EqnBasisTrun}) models can be represented in this way,
for appropriate choices of the operator $\Phi_m$. Recall
representation (\ref{EqnLinopRep}) of $\Phi_m$ in terms of the
functions $\{\phi_j\}_{j=1}^m$.
\begin{itemize}
\item For the time sampling model (\ref{EqnTimeSamp}), we set $\phi_j
= \Kbb(\cdot, t_j)/\sqrt{m}$, so that by the reproducing
property of the kernel, we have $\langle\phi_j,f\rangle_\Hil=
f(t_j)/\sqrt{m}$ for all $f \in\Hil$, and $j = 1, 2, \ldots,
m$. With these choices, the operator $\Phi_m$ maps each
$f \in\Hil$ to the $m$-vector of rescaled samples
\[
\frac{1}{\sqrt{m}} \lleft[\matrix{ f(t_1) & \cdots&
f(t_m) } \rright]^T.
\]
Defining the rescaled noise $\sigma_m=
\frac{\sigma_0}{\sqrt{m}}$ yields an instantiation of
model (\ref{EqnLinObs}) which is equivalent to
time-sampling (\ref{EqnTimeSamp}).
\item For the basis truncation model (\ref{EqnBasisTrun}), we set
$\phi_j = \mu_j \psi_j$ so that the operator $\Phi$ maps each
function $f \in\Hil$ to the vector of basis
coefficients
$ [{\langle\psi_1,f\rangle_{L^2} \enskip\cdots\enskip\langle
\psi_m,f\rangle_{L^2}
} ]^T$.
Setting $\sigma_m= \sigma_0$ then yields
another instantiation of model (\ref{EqnLinObs}), this one
equivalent to basis truncation (\ref{EqnBasisTrun}).
\end{itemize}
A remark on notation before proceeding: in the remainder of the paper,
we use $(\Phi, \sigma)$ as shorthand notation for $(\Phi_m,
\sigma_m)$, since the index $m$ should be implicitly
understood throughout our analysis.

In this paper, we provide and analyze estimators for the
$r$-dimensional eigen-subspace spanned by $\{f^*_j\}$, in both
the sampled domain $\real^m$ and in the functional domain. To be
more specific, for $j = 1,\ldots, r$, define the vectors $z^*_j
:=\Phi f^*_j \in\real^m$, and the subspaces
%
\begin{equation}
\label{eqdefZfsFfs} \mathfrak{Z}^*:=\Span\bigl\{ z^*_1,
\ldots,z^*_r\bigr\} \subset\reals^m \quad\mbox{and}\quad
\mathfrak{F}^*:=\Span\bigl\{ f^*_1,\ldots,f^*_r\bigr\}
\subset\Hil,
\end{equation}
and let $\widehat{\mathfrak{Z}}$ and $\widehat{\mathfrak{F}}$
denote the corresponding estimators. In
order to measure the performance of the estimators, we will use
projection-based distances between subspaces. In particular, let
${P}_{\mathfrak{Z}^*}$ and ${P}_{\widehat{\mathfrak{Z}}}$ be
orthogonal projection
operators into $\mathfrak{Z}^*$ and~$\widehat{\mathfrak{Z}}$,
respectively, considered as
subspaces of $\ell_2^{m}:=
(\reals^m,\|\cdot\|_2)$. Similarly, let ${P}_{\mathfrak
{F}^*}$ and
${P}_{\widehat{\mathfrak{F}}}$ be orthogonal projection operators
into $\mathfrak{F}^*$
and $\widehat{\mathfrak{F}}$, respectively, considered as subspaces of
$(\Hil,\|\cdot\|_{L^2})$. We are interested in bounding the
deviations
%
\begin{equation}
\label{eqprojnormfirstdef}\quad \dist_{\mathrm{HS}}\bigl(\widehat{\mathfrak{Z}},
\mathfrak{Z}^*\bigr):=\nnorm{P}_{\widehat{\mathfrak{Z}}} -
{P}_{\mathfrak{Z}^*}
\nnorm_{\mathrm{HS}}\quad\mbox{and}\quad\dist_{\mathrm{HS}}\bigl(\widehat{
\mathfrak{F}},\mathfrak{F}^*\bigr):=\nnorm{P}_{\widehat{\mathfrak{F}}} -
{P}_{\mathfrak{F}^*} \nnorm_{\mathrm{HS}},
\end{equation}
where $\nnorm\cdot\nnorm_{\mathrm{HS}}$ is the Hilbert--Schmidt norm of an
operator (or matrix).


\subsection{Approximation-theoretic quantities}

One object that plays an important role in our analysis is the matrix
$K:=\Phi\Phi^* \in\real^{m\times m}$. From the
form of the adjoint, it can be seen that $[K]_{ij} =
\langle\phi_i, \phi_j \rangle_\Hil$. For future reference, let us compute
this matrix for the two special cases of linear operators considered
thus far:
\begin{itemize}
\item For the time sampling model (\ref{EqnTimeSamp}), we have $\phi_j
= \Kbb(\cdot, t_j)/\sqrt{m}$ for all $j = 1,\ldots, m$,
and hence
$[K]_{ij} = \frac{1}{m}\langle\Kbb(\cdot, t_i),\Kbb
(\cdot, t_j)\rangle_{\Hil} = \frac{1}{m} \Kbb(t_i, t_j)$,
using the reproducing
property of the kernel.
\item For the basis truncation model (\ref{EqnBasisTrun}), we have
$\phi_j = \mu_j\psi_j$, and hence
%
$[K]_{ij} =
\langle\mu_i \psi_i,\mu_j\psi_j\rangle_{\Hil} = \mu_i\delta_{ij}$.
Thus, in this
special case, we have $K= \diag(\mu_1,\ldots,\break\mu_m)$.
\end{itemize}

In general, the matrix $K$ is a type of Gram matrix, and so is
symmetric and positive semidefinite. We assume\vspace*{1pt} throughout this paper
that the functions $\{\phi_j\}_{j=1}^m$ are linearly independent
in $\Hil$, which implies that $K$ is strictly positive definite.
Consequently, it has a set of eigenvalues which can be ordered as
%
\begin{equation}
\widehat{\mu}_1 \geq\widehat{\mu}_2 \geq\cdots\geq
\widehat{\mu}_m> 0.
\end{equation}
Under this condition, we may use $K$ to define a norm on
$\real^m$ via $\|z\|_K^2:=z^T K^{-1} z$. Moreover, we have
the following interpolation lemma, which is proved in
Appendix~B.1 
of the supplementary material~\cite{supp}:
%
\begin{lemma}
\label{LemInterpolate}
For any $f \in\Hil$, we have $\| \Phi f\|_K \leq
\|f\|_\Hil$, with equality if and only if $f \in\Range(\Phi^*)$.
Moreover, for any $z \in\real^m$, the function $g = \Phi^*
K^{-1} z$ has smallest Hilbert norm of all functions satisfying
$\Phi g = z$, and is the unique function with this property.
\end{lemma}
This lemma is useful in constructing a function-based
estimator, as will be clarified in Section~\ref{SecMest}.

In our analysis of the functional error $\dist_{\mathrm{HS}}(\widehat
{\mathfrak{F}},\mathfrak
{F}^*)$, a
number of approximation-theoretic quantities play an important role.
As a mapping from an infinite-dimensional space $\Hil$ to
$\real^m$, the operator $\Phi$ has a nontrivial nullspace.
Given the observation model (\ref{EqnLinObs}), we receive no
information about any component of a function $f^*$ that lies within
this nullspace. For this reason, we define the width of the nullspace
in the $L^2$-norm, namely the quantity
%
\begin{equation}
\label{EqnDefnPhiNull} N_m(\Phi):=\sup\bigl\{ \|f
\|_{L^2}^2 \mid f \in\Ker(\Phi), \|f\|_\Hil\leq1
\bigr\}.
\end{equation}
In addition, the observation operator $\Phi$ induces a semi-norm on
the space $\Hil$, defined by
%
\begin{equation}
\label{EqnDefnPhinorm} \|f\|_{\Phi}^2:=\| \Phi f
\|_2^2 = \sum_{j=1}^m
[\Phi f]_j^2.
\end{equation}
It is of interest to assess how well this semi-norm approximates the
$L^2$-norm. Accordingly, we define the quantity
%
\begin{equation}
\label{EqnDefnPhiDefect} D_m(\Phi):= \mathop{\sup_{ f \in\Range(\Phi^*)
}}_{ \|f\|_\Hil\leq1}
\bigl| \|f\|_{\Phi}^2 - \|f\|_{L^2}^2 \bigr|,
\end{equation}
which measures the worst-case gap between these two (semi)-norms,
uniformly over the Hilbert ball of radius one, restricted to the
subspace of interest $\Range(\Phi^*)$. Given knowledge of the
linear operator $\Phi$, the quantity $D_m(\Phi)$ can be computed in a
relatively straightforward manner. In particular, recall the
definition of the matrix~$K$, and let us define a second matrix
$\Theta\in{\mathbb{S}^{m}_+}$ with entries $\Theta_{ij}:=
\langle\varphi_i, \varphi_j \rangle_{L^2}$.
%
\begin{lemma}
\label{LemDefect}
We have the equivalence
%
\begin{equation}
D_m(\Phi)= \nnnorml K - K^{-1/2} \Theta K^{-1/2}
\nnnormr_{{2}},
\end{equation}
where \mbox{$\nnorm\cdot\nnorm_{{2}}$} denotes the $\ell_2$-operator norm.
\end{lemma}

See Appendix B.2 
of the supplementary material~\cite{supp} for the proof of this claim.


\section{$M$-estimator and implementation}
\label{SecMest}

With this background in place, we now turn to the description of our
$M$-estimator, as well as practical details associated with its
implementation.

\subsection{$M$-estimator}

We begin with some preliminaries on notation, and our representation
of subspaces. Recall definition (\ref{eqdefZfsFfs}) of
$\mathfrak{Z}^*$ as the $r$-dimen\-sional subspace of $\reals^m$ spanned
by $\{z^*_1,\ldots,z^*_r\}$, where $z^*_j = \Phi f^*_j$.
Our initial goal is to construct an estimate $\widehat{\mathfrak{Z}}$,
itself an $r$-dimensional subspace, of the unknown
subspace~$\mathfrak{Z}^*$.

We represent subspaces by elements of the Stiefel manifold
${V_{r}(\reals^{m})}$, which consists of $m
\times r$ matrices $Z$ with orthonormal columns
\[
{V_{r}\bigl(\reals^{m}\bigr)}:= \bigl\{ Z\in
\reals^{m\times
r} \mid Z^T Z= I_r \bigr\}.
\]
A given matrix $Z$ acts as a representative of the subspace spanned
by its columns, denoted by $\colsp(Z)$. For any $U\in
{V_{r}(\reals^{r})}$, the matrix $ZU$ also belongs to the
Stiefel manifold, and since $\colsp(Z) = \colsp(ZU)$, we may
call $ZU$ a version of $Z$. We let ${P}_Z= ZZ^T
\in\real^{m\times m}$ be the orthogonal projection onto
$\colsp(Z)$. For two matrices $Z_1, Z_2 \in
{V_{r}(\reals^{m})}$, we measure the distance between the
associated subspaces via $\dist_{\mathrm{HS}}(Z_1,Z_2):=
\nnorm{P}_{Z_1}-{P}_{Z_2} \nnorm_{\mathrm{HS}}$, where
$\nnorm\cdot\nnorm_{\mathrm{HS}}$ is the Hilbert--Schmidt (or Frobenius)
matrix norm.

\subsubsection{Subspace-based estimator}

With this notation, we now specify an $M$-estimator for the subspace
$\mathfrak{Z}^*= \myspan\{ z^*_1,\ldots, z^*_r\}$. Let us begin with
some intuition. Given the $n$ samples $\{ \yobs{1},\ldots,
\yobs{n} \}$, let us define the $m\times m$ sample
covariance matrix $\widehat{\Sigma}_{n}:=\frac{1}{{n}}
\sum_{i=1}^{n}\yobs{i} \yobs{i}^T$. Given the observation
model (\ref{EqnLinObs}), a~straightforward computation shows that
%
\begin{equation}
\Exs[\widehat{\Sigma}_{n}] = \sum_{j=1}^rs^2_j
z^*_j \bigl(z^*_j\bigr)^T +
\sigma_m^2 I_{m}.
\end{equation}
Thus, as $n$ becomes large, we expect that the top $r$
eigenvectors of $\widehat{\Sigma}_{n}$ might give a good
approximation to
$\myspan\{z^*_1,\ldots, z^*_r\}$. By the Courant--Fischer
variational representation, these $r$ eigenvectors can be obtained
by maximizing the objective function
\[
\llangle{\widehat{\Sigma}_{n}},{{P}_Z}\rrangle:=\trace
\bigl(\widehat{\Sigma}_{n}Z Z^T\bigr)
\]
over all matrices $Z \in{V_{r}(\reals^{m})}$.

However, this approach fails to take into account the smoothness
constraints that the vectors $z^*_j = \Phi f^*_j$ inherit from the
smoothness of the eigenfunctions $f^*_j$. Since $\|f^*_j\|_\Hil\leq
\rho$ by assumption, Lemma~\ref{LemInterpolate} implies that
\[
\bigl\|z^*_j\bigr\|_K^2 = \bigl(z^*_j
\bigr)^T K^{-1} z^*_j \leq\bigl\|f^*_j
\bigr\|_\Hil^2 \leq\rho^2 \qquad\mbox{for all $j =
1, 2,\ldots, r$.}
\]
Consequently, if we define the matrix $Z^*:= [{ z^*_1
\enskip\cdots\enskip z^*_r
} ] \in\real^{m\times r}$, then it must satisfy
the \textit{trace smoothness condition}
%
\begin{equation}
\label{EqnTraceSmooth} \llanglel{K^{-1}},{Z^* \bigl(Z^*
\bigr)^T}\rrangler= \sum_{j=1}^r
\bigl(z^*_j\bigr)^T K^{-1} z^*_j
\leq r\rho^2.
\end{equation}
This calculation motivates the constraint $\llangle
{K^{-1}},{{P}_Z}\rrangle\leq2 r\rho^2$ in our estimation procedure.

Based on the preceding intuition, we are led to consider the
optimization problem
%
\begin{equation}
\label{EqnDiscreteMest} \widehat{Z}\in\mathop{\arg\max}_{Z\in{V_{r}(\reals
^{m})}} \bigl\{ \llangle{
\widehat{\Sigma}_{{n}}},{{P}_Z}\rrangle\mid
\llanglel{K^{-1}},{{P}_Z}\rrangler\le2 r\rho^2
\bigr\},
\end{equation}
where we recall that ${P}_Z= ZZ^T \in\real^{m
\times
m}$. Given any optimal solution $\widehat{Z}$, we return the subspace
$\widehat{\mathfrak{Z}}= \colsp(\widehat{Z})$ as our estimate of
$\mathfrak{Z}^*$. As
discussed at more
length in Section~\ref{SecImplementation}, it is straightforward to
compute $\widehat{Z}$ in polynomial time. The reader might wonder why we
have included an additional factor of two in this trace smoothness
condition. This slack is actually needed due to the potential
infeasibility of the matrix ${Z^*}$ for
to problem (\ref{EqnDiscreteMest}), which arises since the columns of ${Z^*}$
are not guaranteed to be orthonormal. As shown by our analysis, the
additional slack allows us to find a matrix ${\widetilde{Z}^*}\in
{V_{r}(\reals^{m})}$ that spans the same subspace as ${Z^*}$, and is
also feasible for to problem (\ref{EqnDiscreteMest}). More formally,
we have:
%
\begin{lemma}
\label{LemFeasible}
Under condition (\ref{eqassumpA2}), there exists a matrix
${\widetilde{Z}^*}\in
{V_{r}(\reals^{m})}$ such that
%
\begin{equation}
\Range\bigl({\widetilde{Z}^*}\bigr) = \Range\bigl({Z^*}\bigr)\quad\mbox{and}
\quad\llanglel{K^{-1}},{{\widetilde{Z}^*}\bigl({\widetilde{Z}^*}
\bigr)^T}\rrangler\leq2 r\rho^2.
\end{equation}
\end{lemma}
See Appendix B.3 
of the supplementary material~\cite{supp} for the proof
of this claim.


\subsubsection{\texorpdfstring{The functional estimate $\widehat{\mathfrak{F}}$}{The functional estimate F}}

Having thus obtained an estimate\setcounter{footnote}{1}\footnote{Here,
$\{\widehat{z}_j\}_{j=1}^r\subset\real^m$ is any collection of vectors
that span $\widehat{\mathfrak{Z}}$. As we are ultimately only
interested in the resulting functional ``subspace,'' it does not matter
which particular collection we choose.} $\widehat{\mathfrak{Z}}=
\myspan\{\widehat{z}_1,\ldots, \widehat{z}_r\}$ of $\mathfrak{Z}^*=
\myspan\{ z^*_1,\ldots, z^*_r\}$, we now need to construct a
$r$-dimen\-sional subspace $\widehat{\mathfrak{F}}$ of the Hilbert space
to be used as an estimate of $\mathfrak{F}^*= \myspan\{ f^*_1,\ldots,
f^*_r\}$. We do so using the interpolation suggested by
Lemma~\ref{LemInterpolate}. For each $j = 1,\ldots, r$, let us define
the function
%
\begin{equation}
\label{EqnDefnFhat} \widehat{f}_j:=\Phi^* K^{-1}
\widehat{z}_j = \sum_{i=1}^m
\bigl(K^{-1} \widehat{z}_j\bigr)_i
\phi_i.
\end{equation}
Since $K= \Phi\Phi^*$ by definition, this construction ensures
that $\Phi\widehat{f}_j = \widehat{z}_j$. Moreover, Lemma \ref
{LemInterpolate}
guarantees that $\widehat{f}_j$ has the minimal Hilbert norm (and
hence is
smoothest in a certain sense) over all functions that have this
property. Finally, since $\Phi$ is assumed to be surjective
(equivalently, $K$ assumed invertible), $\Phi^* K^{-1}$ maps
linearly independent vectors to linearly independent functions, and
hence preserves dimension. Consequently, the space $\widehat{\mathfrak
{F}}:=\myspan
\{
\widehat{f}_1,\ldots, \widehat{f}_r\}$ is an $r$-dimensional
subspace of
$\Hil$ that we take as our estimate of $\mathfrak{F}^*$.



\subsection{Implementation details}
\label{SecImplementation}

In this section, we consider some practical aspects of implementing
the $M$-estimator, and present some simulations to illustrate its
qualitative properties. We begin by observing that once the subspace
vectors $\{\widehat{z}_j\}_{j=1}^r$ have been computed, then it is
straightforward to compute the function estimates $\{ \widehat{f}_j
\}_{j=1}^r$, as weighted combinations of the functions
$\{\phi_j\}_{j=1}^m$. Accordingly, we focus our attention on
solving problem~(\ref{EqnDiscreteMest}).

On the surface, problem (\ref{EqnDiscreteMest}) might appear
nonconvex, due to the Stiefel manifold constraint. However, it
can be reformulated as a semidefinite program (SDP), a well-known
class of convex programs, as clarified in the following:
%
\begin{lemma}
\label{LemImplement}
Problem (\ref{EqnDiscreteMest}) is equivalent to solving the SDP
%
\begin{eqnarray}
\label{EqnSDP} \widehat{X}\in\mathop{\arg\max}_{X \succeq0 } \llangle{\widehat
{\Sigma
}_n},{X}\rrangle\hspace*{45pt}\nonumber\\[-8pt]\\[-8pt]
&&\eqntext{\mbox{such that $\nnorm X
\nnorm_{{2}} \leq1$, $\trace(X) = r$, and $\llanglel{K^{-1}},{X}
\rrangler\leq2 r\rho^2$}}
\end{eqnarray}
for which there always exists an optimal rank $r$ solution.
Moreover, by Lagrangian duality, for some $\beta> 0$, the problem is
equivalent to
%
\begin{equation}
\label{eqsdpregver}\qquad \widehat{X}\in\mathop{\arg\max}_{X \succeq0 } \llanglel
{\widehat{
\Sigma}_n- \beta K^{-1}},{X}\rrangler\qquad\mbox{such that
$\nnorm X \nnorm_{{2}} \leq1$ and $\trace(X) = r$},
\end{equation}
which can be solved by an eigen decomposition of $\widehat{\Sigma}_n-
\beta K^{-1}$.
\end{lemma}

As a consequence, for a given Lagrange multiplier $\beta$,
the regularized form of the estimator can be solved with the cost of
solving an eigenvalue problem. For a given constraint $2 r
\rho^2$, the appropriate value of $\beta$ can be found by a
path-tracing algorithm, or a simple dyadic splitting approach.

In practice where the radius $\rho$ is not known, one could
use cross-validation to set a proper value for the Lagrange multiplier
$\beta$. A possibly simpler approach is to evaluate
$\llangle{K^{-1}},{X}\rrangle$ for the optimal $X$ on a grid of
$\beta$ and
choose a value around which $\llangle{K^{-1}},{X}\rrangle$ is least
variable. As for the choice of the number of components~$r$, a\vadjust{\goodbreak}
standard approach for choosing it would be to compute the estimator
for different choices, and plot the residual sum of eigenvalues of the
sample covariance matrix. As in ordinary PCA, an elbow in such a plot
indicates a proper trade-off between the number of components to
keep and the amount of variation explained.

\begin{figure}

\includegraphics{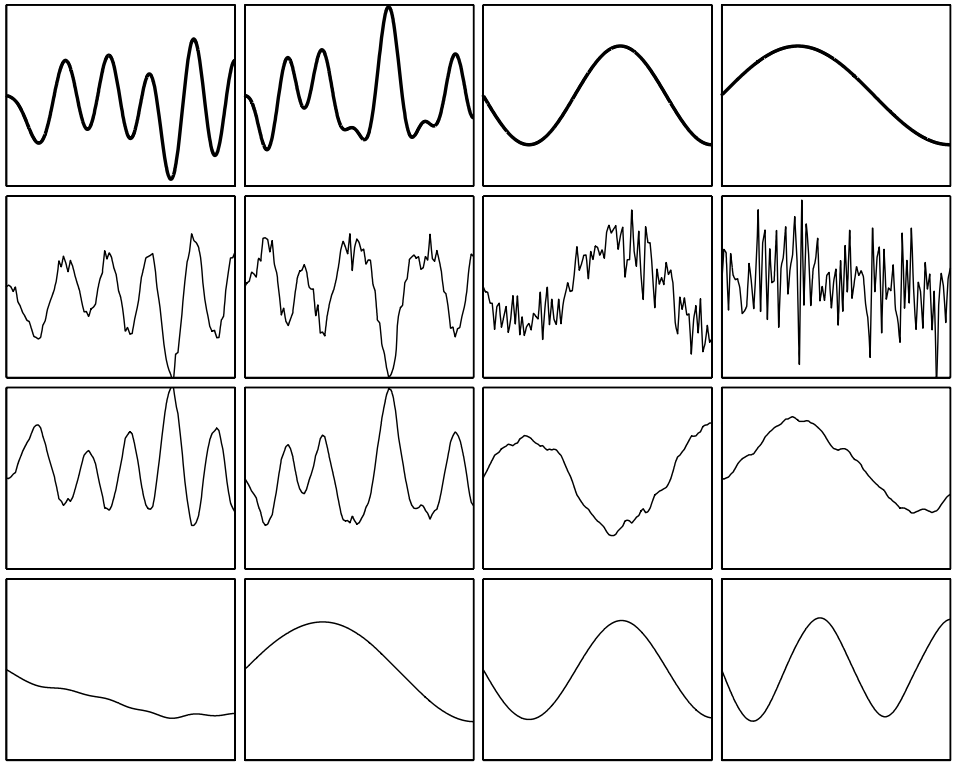}

\caption{Regularized PCA for time sampling in first-order Sobolev
RKHS. Top row shows, from left to right, plots of the $r=4$
``true'' principal components $f^*_{1},\ldots,f^*_{4}$ with
signal-to-noise ratios $s_{1} = 1, s_{2} = 0.5, s_{3} = 0.25$
and $s_{4} = 0.125$, respectively. The
number of statistical and functional samples are $n= 75$ and
$m= 100$. Subsequent rows show the corresponding estimators
$\widehat{f}_{1},\ldots, \widehat{f}_{4}$ obtained by applying the regularized
form (\protect\ref{eqsdpregver}).}
\label{fig1}
\end{figure}

In order to illustrate the estimator, we consider
the time sampling model (\ref{EqnTimeSamp}), with uniformly spaced samples,
in the context of a first-order Sobolev RKHS [with kernel function
$\Kbb(s,t) = \min(s,t)$]. The
parameters of the model are taken to be $r= 4$,
$(s_1,s_2,s_3,s_4) = (1, 0.5, 0.25, 0.125)$,
$\sigma_0 = 1$, $m= 100$ and $n= 75$. The regularized form
(\ref{eqsdpregver}) of the estimator is applied, and the results are
shown in Figure~\ref{fig1}.
The top row corresponds to the four ``true'' signals $\{f^*_j\}$, the
leftmost being $f^*_1$ (i.e., having the highest signal-to-noise
ratio)\vspace*{1pt} and the rightmost $f^*_4$. The subsequent rows show
the corresponding estimates $\{\widehat{f}_j\}$, obtained using
different values of $\beta$. The second, third and fourth rows
correspond to $\beta =0$, $\beta= 0.0052$ and $\beta= 0.83$.

One observes that without regularization ($\beta= 0$), the
estimates for the two weakest signals ($f^*_3$ and $f^*_4$) are poor. The
case $\beta=0.0052$ is roughly the one which achieves the minimum
for the dual problem. One observes that the quality of the estimates
of the signals, and in particular the weakest ones, are considerably
improved. The optimal (oracle) value of $\beta$, that is, the one
which achieves the minimum error between $\{f^*_j\}$ and $\{\widehat
{f}_j\}$,
is $\beta= 0.0075$ in this problem. The corresponding estimates
are qualitatively similar to those of $\beta=0.0052$ and are not
shown.

The case $\beta= 0.83$ shows the effect of over-regularization. It
produces very smooth signals, and although it fails to reveal $f^*_1$
and $f^*_2$, it reveals highly accurate versions of $f^*_3$ and
$f^*_4$. It is also interesting to note that the smoothest signal,
$f^*_4$, now occupies the position of the second (estimated) principal
component. That is, the regularized PCA sees an effective
signal-to-noise ratio which is influenced by smoothness.
This suggests a rather practical appeal of the method in
revealing smooth signals embedded in noise. One can vary $\beta$
from zero upward, and if some patterns seem to be present for a wide
range of $\beta$ (and getting smoother as $\beta$ is
increased), one might suspect that they are indeed present in data but
masked by noise.


\section{Main results}
\label{SecMain}

We now turn to the statistical analysis of our estimators, in particular
deriving high-probability upper bounds on the error of the subspace-based
estimate $\widehat{\mathfrak{Z}}$, and the functional estimate
$\widehat{\mathfrak{F}}$. In both cases,
we begin by stating general theorems that apply to arbitrary linear
operators $\Phi$---Theorems~\ref{THMSUBSPACEERR} and~\ref{THMFUNCERR},
respectively---and then derive a number of corollaries for particular
instantiations of the observation operator.

\subsection{\texorpdfstring{Subspace-based estimation rates (for $\widehat{\mathfrak{Z}}$)}{Subspace-based estimation rates (for Z)}}
\label{secconsistZfh}

We begin\vspace*{1pt} by stating high-probability upper bounds on the
error $\dist_{\mathrm{HS}}(\widehat{\mathfrak{Z}}, \mathfrak{Z}^*)$ of
the subspace-based estimates. Our rates are stated in terms of a
function that involves the eigenvalues of the matrix $K= \Phi\Phi^*
\in\real^m$, ordered as $\widehat{\mu}_1 \ge\widehat{\mu}_2
\geq\cdots\geq\widehat{\mu}_m> 0$. Consider the function $\SimpleF\dvtx
\real_+ \rightarrow\real_+$ given by
%
\begin{equation}
\label{EqnGrowthFuncDef} \SimpleF(t):= \Biggl[ \sum
_{j=1}^m\min\bigl\{t^2, r
\rho^2 \widehat{\mu}_j \bigr\} \Biggr]^{1/2}.
\end{equation}
As will be clarified in our proofs, this function provides a measure
of the statistical complexity of the function
class
\[
\Range\bigl(\Phi^*\bigr) = \Biggl\{f \in\Hil\Bigm| f = \sum
_{j=1}^ma_j \phi_j
\mbox{ for some $a \in\real^m$} \Biggr\}.
\]

We require a few regularity assumptions. Define the quantity
%
\begin{equation}
\label{EqnDefnPhiAngle} C_m\bigl(f^*\bigr):=\max_{1 \le i,j \le r} \bigl|
\bigl\langle f^*_i,f^*_j\bigr\rangle_{\Phi} -
\delta_{ij} \bigr| = \max_{1 \le
i,j \le r} \bigl| \bigl\langle z^*_i,z^*_j
\bigr\rangle_{\real^m} - \delta_{ij} \bigr|,
\end{equation}
which measures\vspace*{1pt} the departure from orthonormality of the vectors $z^*_j
:=\Phi f^*_j$ in~$\real^m$. A straightforward argument
using a polarization identity shows that $C_m(f^*)$ is upper bounded (up
to a constant factor) by the uniform quantity $D_m(\Phi)$, as defined in
equation (\ref{EqnDefnPhiDefect}). Recall that the random functions
are generated according to the model $\xobs{i} = \sum_{j=1}^r
s_j \beta_{ij} f^*_j$, where the signal strengths are ordered as
$1 = s_1 \geq s_2 \geq\cdots\geq s_r> 0$,
and that $\sigma_m$ denotes the noise standard deviation in the
observation model (\ref{EqnLinObs}).

In terms of these quantities, we require the following assumptions:
%
\begin{subequation}
%
\begin{eqnarray}
&&
\mbox{(A1)}\quad  \frac{s_r^2}{s_1^2} \geq\frac{1}{2}\quad\mbox{and}\quad
\sigma_0^2 := \sup_{m} \sigma_m^2
\leq\kappa s_1^2,
\\
\label{eqassumpA2}
&&
\mbox{(A2)}\quad  C_m\bigl(f^*\bigr) \leq
\frac{1}{2r} \quad\mbox{and}
\\
&&
\mbox{(A3)}\quad \frac{\sigma_m}{\sqrt{n}} \SimpleF(t) \leq\sqrt{\kappa} t
\qquad\mbox{for the same constant $\kappa$ as in (A1),}
\\
\label{eqassumpA4}
&&
\mbox{(A4)}\quad r \le \min\biggl\{ \frac{m}{2},
\frac{n}{4}, \kappa\frac{\sqrt{n}}{\sigma_m
} \biggr\}.
\end{eqnarray}
\end{subequation}
%

\begin{Remarks*}
The first part of condition (A1) is to prevent the ratio $s_r
/ s_1$ from going to zero as the pair $(m, n)$
increases, where the constant $1/2$ is chosen for convenience. Such a
lower bound is necessary for consistent estimation of the
eigen-subspace corresponding to $\{s_1,\ldots,s_r\}$. The
second part of condition (A1), involving the constant $\kappa$,
provides a lower bound on the signal-to-noise ratio $s_r/
\sigma_m$. Condition (A2) is required\vspace*{1pt} to prevent degeneracy among
the vectors $z^*_j = \Phi f^*_j$ obtained by mapping the unknown
eigenfunctions to the observation space $\real^m$. [In the ideal
setting, we would have $C_m(f^*)= 0$, but our analysis shows that the
upper bound in (A2) is sufficient.] Condition (A3) is required so
that the critical tolerance $\epscrit$ specified below is
well-defined; as will be clarified, it is always satisfied for the
time-sampling model, and holds for the basis truncation model whenever
$n\geq m$. Condition (A4) is easily satisfied, since the
right-hand side of (\ref{eqassumpA4}) goes to infinity while we
usually take
$r$ to be fixed. Our results, however, are still valid if $r$
grows slowly with $m$ and $n$ subject
to (\ref{eqassumpA4}).
\end{Remarks*}
%
\begin{theorem}
\label{THMSUBSPACEERR}
Under conditions \textup{(A1)--(A3)} for a sufficiently small constant
$\kappa$, let $\epscrit$ be the smallest positive number satisfying
the inequality
%
\begin{equation}
\label{EqnDefnEpscrit} \frac{\sigma_m}{\sqrt{{n}}} r^{3/2} \SimpleF
(\epsilon)
\leq\kappa\epsilon^2.
\end{equation}
Then there are universal positive constants $(c_0,
c_1, c_2)$ such that
%
\begin{equation}
\mprob\bigl[ \dist^2_{\mathrm{HS}}\bigl(\widehat{\mathfrak{Z}},
\mathfrak{Z}^*\bigr) \leq c_0 \epscrit^2 \bigr] \geq1 -
\varphi(n, \epscrit),
\end{equation}
where $\varphi(n, \epscrit):=c_1 \{
r^2 \exp( - c_2 r^{-3} \frac{{n}
}{\sigma_m^2} (\epscrit\wedge
\epscrit^2) ) + r\exp( - \frac{{n}}{64}) \}$.
\end{theorem}


We note that Theorem~\ref{THMSUBSPACEERR} is a general result,
applying to an arbitrary bounded linear operator $\Phi$. However,
we can obtain a number of concrete results by making specific choices
of this sampling operator, as we explore in the following sections.

\subsubsection{Consequences for time-sampling}

Let us begin with the time-sam\-pling model (\ref{EqnTimeSamp}), in
which we observe the sampled functions
\[
y_i = \lleft[\matrix{ x_i(t_1) &
x_i(t_2) & \cdots& x_i(t_m) }
\rright]^T + \sigma_0 w_i\qquad
\mbox{for $i = 1, 2,\ldots, m$.}
\]
As noted earlier, this set-up can be modeled in our general
setting (\ref{EqnLinObs}) with $\phi_j = \Kbb(\cdot,
t_j)/\sqrt{m}$ and $\sigma_m= \sigma_0/\sqrt{m}$.

In this case,\vspace*{1pt} by the reproducing property of the RKHS, the matrix $K
= \Phi\Phi^*$ has entries of the form $K_{ij} =
\langle\phi_i,\phi_j\rangle_{\Hil} = \frac{\Kbb(t_i, t_j)}{m}$.
Letting $\widehat{\mu}_1 \geq\widehat{\mu}_2 \geq\cdots\geq
\widehat{\mu}_m>
0$ denote its ordered eigenvalues, we say that the kernel matrix $K$
has polynomial-decay with parameter $\alpha> 1/2$ if there is a
constant $c$ such that $\widehat{\mu}_j \leq c j^{-2\alpha}$ for
all $j
= 1, 2,\ldots, m$. Since the kernel matrix $K$ represents a
discretized approximation of the kernel integral operator defined by
$\Kbb$, this type of polynomial decay is to be expected whenever the
kernel operator has polynomial-$\alpha$ decaying eigenvalues. For
example, the usual spline kernels that define Sobolev spaces have this
type of polynomial decay~\cite{Gu02}. In Appendix A 
of the supplementary material~\cite{supp}, we verify this property
explicitly for the kernel $\Kbb(s,t) = \min\{s, t\}$ that defines the
Sobolev class with smoothness $\alpha= 1$.

For any such kernel, we have the following consequence of
Theorem~\ref{THMSUBSPACEERR}:
%
\begin{corollary}[(Achievable rates for time-sampling)]
\label{CorTimeEmp} Consider the case of a time-sampling operator
$\Phi$. In addition to conditions \textup{(A1)} and \textup{(A2)},
suppose that the kernel matrix $K$ has polynomial-decay with parameter
$\alpha> 1/2$. Then we have
%
\begin{equation}
\label{EqnTimeEmp}\quad \mprob\biggl[ \dist^2_{\mathrm{HS}}\bigl(
\widehat{\mathfrak{Z}}, \mathfrak{Z}^*\bigr) \leq c_0 \min\biggl\{
\biggl( \frac{\AnnoyCon \sigma_0^2}{mn} \biggr)^{{
{2 \alpha}/({2 \alpha+1})}}, r^3
\frac{\sigma_0^2}{n} \biggr\} \biggr] \geq1 - \varphi(n, m),
\end{equation}
where $\AnnoyCon:=r^{3+{1}/({2\alpha})}
\rho^{{1}/{\alpha}}$, and $\varphi(n, m):=
c_1 \{ \exp(- c_2 \{ (
r^{-2}
\rho^2 mn )^{{1}/({2 \alpha+1})} \wedge m
\} ) + \exp( -n/64) \}$.
\end{corollary}

\begin{Remarks*}
(a) Disregarding constant pre-factors not depending on the pair
$(m, n)$, Corollary~\ref{CorTimeEmp} guarantees\vspace*{1pt} that solving
problem (\ref{EqnDiscreteMest}) returns a subspace estimate $\widehat
{\mathfrak{Z}}$
such that
\[
\dist^2_{\mathrm{HS}}\bigl(\widehat{\mathfrak{Z}}, \mathfrak{Z}^*
\bigr) \precsim\min\bigl\{(mn )^{-{{2
\alpha}/({2 \alpha+1})}}, n^{-1} \bigr\}
\]
with high probability as $(m,
n)$ increase. Depending on the scaling of the number of
time samples $m$ relative to the number of functional samples
$n$, either term\vadjust{\goodbreak} in this upper bound can be the smallest (and
hence active) one. For instance, it can be verified that whenever
$m\geq
n^{{1}/({2 \alpha})}$, then the first term is smallest, so
that we achieve the rate $\dist^2_{\mathrm{HS}}(\widehat{\mathfrak
{Z}}, \mathfrak{Z}^*)
\precsim
(mn)^{-{{2 \alpha}/({2 \alpha+1})}}$. The appearance of the term $(m
n)^{-{{2 \alpha}/({2 \alpha+1})}}$ is quite natural, as it
corresponds to the
minimax rate of a nonparametric regression problem with smoothness
$\alpha$, based on $m$ samples each of variance
$n^{-1}$.
Later, in Section~\ref{SecLower}, we provide results guaranteeing that
this scaling is minimax optimal under reasonable conditions on the
choice of sample points; in particular, see
Theorem~\ref{THMMINIMAXEMP}(a).

(b) To be clear, although bound (\ref{EqnTimeEmp}) allows for the
possibility that the error is of order \textit{lower than
$n^{-1}$}, we note that the probability with which the
guarantee holds includes a term of the order $\exp(-n/64)$.
Consequently, in terms of expected error, we cannot guarantee a rate
faster than $n^{-1}$.
\end{Remarks*}
\begin{pf*}{Proof of Corollary~\ref{CorTimeEmp}}
We need to bound the critical value $\epscrit$ defined in the theorem
statement (\ref{EqnDefnEpscrit}). Define the function
$\mathcal{G}^2(t):=\break \sum_{j=1}^m\min\{ \widehat{\mu}_j, t^2
\}$, and note that $\SimpleF(t) = \sqrt{ r} \rho
\mathcal{G}(\frac{t}{\sqrt{r} \rho})$ by
construction.
Under
the assumption of polynomial-$\alpha$ eigendecay, we have
\[
\mathcal{G}^2(t) \leq\int_0^\infty
\min\bigl\{ c x^{-2\alpha}, t^2 \bigr\} \,dx,
\]
and some algebra then shows that $\mathcal{G}(t) \precsim t^{1 -
1/(2\alpha)}$.
Disregarding constant factors, an upper bound on the critical
$\epscrit$ can be obtained by solving the equation
\[
\epsilon^2 = \frac{\sigma_m}{\sqrt{{n}}} r^{3/2} \sqrt{r} \rho
\biggl(\frac{\epsilon}{\sqrt{r}\rho} \biggr)^{1-1/(2\alpha)}.
\]
Doing so yields the upper bound $\epsilon^2 \precsim[
\frac{\sigma_m^2}{{n}} r^{3} (\sqrt{r}
\rho)^{{1}/{\alpha}} ]^{{{2 \alpha}/({2 \alpha
+1})}}$. Otherwise, we also
have the trivial upper bound $\SimpleF(t) \le\sqrt{m} t$, which
yields the alternative upper bound $\varepsilon_{m,{n}} \precsim
(\frac{m \sigma_m^2}{{n}}r^3 )^{1/2}$.
Recalling that $\sigma_m= \sigma_0/\sqrt{m}$ and combining the
pieces yields the claim. Notice that this last (trivial) bound on
$\SimpleF(t)$
implies that condition (A3) is always satisfied for the time-sampling
model.
\end{pf*}
%

\subsubsection{Consequences for basis truncation}

We now turn to some consequences for the basis truncation
model (\ref{EqnBasisTrun}).
%
\begin{corollary}[(Achievable rates for basis truncation)]
\label{CorTrunEmp}
Consider a basis truncation operator $\Phi$ in a Hilbert space with
polynomial-$\alpha$ decay. Under conditions \textup{(A1), (A2)} and $m\le
n$, we have
%
\begin{equation}
\mprob\biggl[ \dist^2_{\mathrm{HS}}\bigl(\widehat{\mathfrak{Z}},
\mathfrak{Z}^*\bigr) \leq c_0 \biggl(\frac{\AnnoyCon \sigma_0^2}{n}
\biggr)^{{{2 \alpha}/({2
\alpha+1})}} \biggr] \geq1 - \varphi(n, m),
\end{equation}
where $\AnnoyCon:=r^{3+{1}/({2\alpha})}
\rho^{{1}/{\alpha}}$, and $\varphi(n, m):=
c_1 \{ \exp(- c_2 ( r^{-2} \rho^2
n )^{{1}/({2 \alpha+1})} ) + \exp( -n/64)
\}$.\vadjust{\goodbreak}
\end{corollary}
\begin{pf}%
We note that as long as $m\le n$, condition (A3) is
satisfied, since $\frac{\sigma_m}{\sqrt{n}} \SimpleF(t)
\le\sigma_0\sqrt{\frac{m}{n}} t \le\sigma_0 t$. The rest
of the proof follows that of Corollary~\ref{CorTimeEmp}, noting that
in the last step we have $\sigma_m= \sigma_0$ for the basis
truncation model.\vspace*{-2pt}
\end{pf}
%



\subsection{\texorpdfstring{Function-based estimation rates (for $\widehat{\mathfrak{F}}$)}{Function-based estimation rates (for F)}}
\label{secconsistFfh}

As mentioned earlier, given the consistency of $\widehat{\mathfrak
{Z}}$, the consistency
of $\widehat{\mathfrak{F}}$ is closely related to approximation
properties of the
semi-norm ${\|\cdot\|_{\Phi}}$ induced by $\Phi$, and in
particular how
closely it approximates the $L^2$-norm. These approximation-theoretic
properties are captured in part by the nullspace width $N_m(\Phi)$ and
defect $D_m(\Phi)$ defined earlier in equations (\ref{EqnDefnPhiNull})
and (\ref{EqnDefnPhiDefect}), respectively. In addition to these
previously defined quantities, we require bounds on the following
global quantity:
%
\begin{equation}
\label{eqdefrpepsnu} \rp(\epsilon; \nu):=\sup\bigl\{ \|f
\|_{L^2}^2 \mid\|f\|_\Hil^2 \leq
\nu^2, \|f\|_{\Phi}^2 \leq\epsilon^2
\bigr\}.
\end{equation}
A general upper bound on this quantity is of the form
%
\begin{equation}
\label{equpboundrpremain} \rp(\epsilon; \nu) \le c_1
\epsilon^2 + \nu^2 S_m(\Phi).
\end{equation}
In fact, it is not hard to show that such a bound exists with $c_1 =2$
and $S_m(\Phi)= 2(D_m(\Phi)+ N_m(\Phi))$ using the decomposition
$\Hil=
\Range(\Phi^*) \oplus\Null(\Phi)$. However, this bound is not
sharp. Instead, one can show that in most cases of interest, the term
$S_m(\Phi)$ is of the order of $N_m(\Phi)$.

There are a variety of conditions that ensure that $S_m(\Phi)$ has this
scaling; we refer the reader to the paper~\cite{AmiWaiApprox} for a
general approach. Here we provide a simple sufficient condition,
namely,
%
\begin{equation}
\mbox{(B1)}\quad \Theta\preceq c_0 K^2
\end{equation}
for a positive constant $c_0$.\vspace*{-2pt}
%
\begin{lemma}\label{lemREMAIN}
Under \textup{(B1)}, bound (\ref{equpboundrpremain}) holds with
$c_1 = 2c_0$ and $S_m(\Phi)= 2N_m(\Phi)$.\vspace*{-2pt}
\end{lemma}

See Appendix B.4 
of the supplementary material~\cite{supp} for
the proof of this claim. In the sequel, we show that the first-order
Sobolev RKHS satisfies condition (B1).\vspace*{-2pt}


%
\begin{theorem}
\label{THMFUNCERR}
Suppose that condition \textup{(A1)} holds, and the approximation-theoretic
quantities satisfy the bounds $D_m(\Phi)\leq\frac1{4r\rho^2}\le
1$ and $N_m(\Phi)\le1$. Then there is a constant ${\kappa'_{r,
\rho}}$ such
that
%
\begin{equation}
\label{eqprojdistfuncdist} \dist_{\mathrm{HS}}^2\bigl(\widehat{
\mathfrak{F}},\mathfrak{F}^*\bigr) \leq{\kappa'_{r, \rho
}}
\bigl\{ \epscrit^2 + S_m(\Phi)+ \bigl[D_m(
\Phi)\bigr]^2 \bigr\}
\end{equation}
with the same probability as in Theorem~\ref{THMSUBSPACEERR}.\vspace*{-2pt}
\end{theorem}
%

As with Theorem~\ref{THMSUBSPACEERR}, this is a generally applicable
result, stated in abstract form. By specializing it to different
sampling models, we can obtain concrete rates, as illustrated in the
following sections.\vadjust{\goodbreak}

\subsubsection{Consequences for time-sampling}

We begin by returning to the case of the time sampling
model (\ref{EqnTimeSamp}), where $\phi_j = \Kbb(\cdot,
t_j)/\sqrt{m}$. In this case, condition (B1) needs to be
verified by some calculations. For instance, as shown in
Appendix A 
of the supplementary material~\cite{supp}, in the case of the Sobolev kernel with
smoothness $\alpha= 1$ [namely, $\Kbb(s,t) = \min\{s, t\}$], we are
guaranteed that (B1) holds with $c_0 = 1$, whenever the samples
$\{t_j\}$ are chosen uniformly over $[0,1]$; hence, by
Lemma~\ref{lemREMAIN}, $S_m(\Phi)= 2N_m(\Phi)$.
Moreover, in the case of uniform sampling, we expect that the
nullspace width $N_m(\Phi)$ is upper bounded by $\mu_{m+1}$, and so
will be proportional to $m^{-2\alpha}$ in the case of a kernel
operator with polynomial-$\alpha$ decay. This is verified
in~\cite{AmiWaiApprox} (up to a logarithmic factor) for the case of
the first-order Sobolev kernel. In Appendix A 
of the supplementary material~\cite{supp}, we also show that, for this kernel, $[D_m(\Phi)]^2$ is of
the order $m^{-2\alpha}$, that is, of the same order as
$N_m(\Phi)$.

%
\begin{corollary}
\label{CorTimeFunc}
Consider the basis truncation model (\ref{EqnBasisTrun}) with
uniformly spaced samples, and assume condition \textup{(B1)} holds and that
$N_m(\Phi)+\break [D_m(\Phi)]^2 \precsim m^{-2\alpha}$. Then the
$M$-estimator returns a subspace estimate $\widehat{\mathfrak{F}}$
such that
%
\begin{equation}
\dist_{\mathrm{HS}}^2\bigl(\widehat{\mathfrak{F}},\mathfrak{F}^*
\bigr) \leq{\kappa'_{r, \rho
}} \biggl\{ \min\biggl\{ \biggl(
\frac{\sigma_0^2}{nm} \biggr)^{{{2 \alpha}/({2 \alpha
+1})}}, \frac{\sigma_0^2}{n} \biggr\} +
\frac{1}{m^{2 \alpha}} \biggr\}
\end{equation}
with the same probability as in Corollary~\ref{CorTimeEmp}.
\end{corollary}

In this case, there is an interesting trade-off between the bias or
approximation error which is of order $m^{-2 \alpha}$ and
the estimation error. An interesting transition occurs at the point
when $m\succsim n^{{1}/({2 \alpha})}$, at which:
\begin{itemize}
\item the bias term $m^{-2 \alpha}$ becomes of the order
$n^{-1}$, so that it is no longer dominant, and
\item for the two terms in the estimation error, we have the ordering
\[
(mn)^{-{{2 \alpha}/({2 \alpha+1})}} \leq\bigl(n^{1 + {1}/({2
\alpha})} \bigr)^{-{{2 \alpha}/({2 \alpha+1})}} =
n^{-1}.
\]
\end{itemize}
Consequently, we conclude that the scaling $m=
n^{{1}/({2 \alpha})}$ is the minimal number of samples such
that we achieve an overall bound of the order $n^{-1}$ in the
time-sampling model. In Section~\ref{SecLower}, we will see that
these rates are minimax-optimal.

\subsubsection{Consequences for basis truncation}

For the basis truncation operator~$\Phi$, we have $\Theta= K^2
= \diag(\mu_1^2,\ldots,\mu_m^2)$ so that condition (B1) is
satisfied trivially with $c_0 = 1$. Moreover, Lemma~\ref{LemDefect}
implies $D_m(\Phi)= 0$. In addition, a function $f
= \sum_{j=1}^\infty\sqrt{\mu_j} a_j \psi_j$ satisfies $\Phi f = 0$
if and only if $a_1 = a_2 = \cdots= a_m= 0$, so that
\[
N_m(\Phi)= \sup\bigl\{ \|f\|_{L^2}^2 \mid
\|f\|_\Hil\leq1, \Phi f = 0 \bigr\} = \mu_{m+1}.
\]
%
Consequently, we obtain the following corollary of
Theorem~\ref{THMFUNCERR}:
%
\begin{corollary}
\label{CorBasisTrunFunc}
Consider the basis truncation model (\ref{EqnBasisTrun}) with a kernel
operator that has polynomial-$\alpha$ decaying eigenvalues. Then
the $M$-estimator returns a function subspace estimate $\widehat
{\mathfrak{F}}$ such
that
%
\begin{equation}
\dist_{\mathrm{HS}}^2\bigl(\widehat{\mathfrak{F}},\mathfrak{F}^*
\bigr) \leq{\kappa'_{r, \rho
}} \biggl\{ \biggl(
\frac{\sigma_0^2}{n} \biggr)^{{{2 \alpha}/({2 \alpha+1})}} + \frac
{1}{m^{2
\alpha}} \biggr\}
\end{equation}
with the same probability as in Corollary~\ref{CorTrunEmp}.
\end{corollary}
By comparison to Corollary~\ref{CorTimeFunc}, we see that the
trade-offs between $(m, n)$ are very different for basis
truncation. In particular, there is \textit{no interaction} between the
number of functional samples $m$ and the number of statistical
samples $n$. Increasing $m$ only reduces the approximation
error, whereas increasing $n$ only reduces the estimation error.
Moreover, in contrast to the time sampling model of
Corollary~\ref{CorTimeFunc}, it is impossible to achieve the fast rate
$n^{-1}$, regardless of how we choose the pair $(m,
n)$. In Section~\ref{SecLower}, we will also see that the rates
given in Corollary~\ref{CorBasisTrunFunc} are minimax optimal.


\subsection{Lower bounds}
\label{SecLower}

We now turn to lower bounds on the minimax risk, demonstrating the
sharpness of our achievable results in terms of their scaling with
$(m, n)$. In order to do so, it suffices to consider the
simple model with a single functional component $f^*\in
\Ball_\Hil(1)$, so that we observe $y_i = \beta_{i1}
\Phi_m(f^*) + \sigma_mw_i$ for $i = 1, 2,\ldots,
n$, where $\beta_{i1} \sim N(0,1)$ are i.i.d. standard normal
variates. The minimax risk over the unit ball of the function space
$\Hil$ in the $\Phi$-norm is given by
%
\begin{equation}
\label{EqnDefnMiniMaxPhi} \MiniMax^\Hil\bigl(\| \cdot\|_{\Phi}\bigr):=
\inf_{\widetilde{f}} \sup_{f^*
\in\Ball_\Hil(1)} \Exs\bigl\| \widetilde{f}- f^*
\bigr\|_{\Phi}^2, 
\end{equation}
where the\vspace*{1pt} function $f^*$ ranges over the unit ball
$\Ball_\Hil(1) = \{f \in\Hil \mid \|f\|_\Hil\leq1 \}$ of some Hilbert
space, and $\widetilde{f}$ ranges over measurable functions of the data
matrix $(y_1, y_2,\ldots, y_n) \in\real^{m\times n}$.



%
\begin{theorem}[(Lower bounds for $\|\widetilde{f}- f^*\|_\Phi$)]
\label{THMMINIMAXEMP}
Suppose that the kernel matrix $K$ has eigenvalues with
polynomial-$\alpha$ decay and \textup{(A1)} holds.
\begin{longlist}[(a)]
\item[(a)] For the time-sampling model, we have
%
\begin{equation}
\MiniMax^\Hil\bigl( \| \cdot\|_\Phi\bigr) \geq C\min\biggl\{
\biggl(\frac{\sigma^2_0}{mn} \biggr)^{{{2 \alpha}/({2
\alpha+1})}
}, \frac{\sigma_0^2}{n} \biggr\}.
\end{equation}

\item[(b)] For the frequency-truncation model, with $m\ge(c_0
n)^{{1}/({2\alpha+1})}$, we have
%
\begin{equation}
\MiniMax^\Hil\bigl( \|\cdot\|_\Phi\bigr) \geq C \biggl(
\frac{\sigma_0^2}{ n} \biggr)^{{{2 \alpha}/({2 \alpha+1})}}. 
\end{equation}
\end{longlist}
\end{theorem}


Note that part (a) of Theorem~\ref{THMMINIMAXEMP} shows that the rates
obtained in Corollary~\ref{CorTimeFunc} for the case of time-sampling
are minimax optimal. Similarly,\vadjust{\goodbreak} comparing part (b) of the theorem to
Corollary~\ref{CorBasisTrunFunc}, we conclude that the rates obtained
for frequency truncation model are minimax optimal for $n\in
[m,\break c_1 m^{2\alpha+ 1}]$. As will become clear momentarily
(as a consequence of our next theorem), the case $n> c_1
m^{2\alpha+ 1}$ is not of practical interest.


We now turn to lower bounds on the minimax risk in the
${\|\cdot\|_{L^2}}$ norm---namely
%
\begin{equation}
\label{EqnDefnMiniMaxLtwo} \MiniMax^\Hil\bigl(\| \cdot\|_{L^2}\bigr):=
\inf_{\widetilde{f}} \sup_{f^*\in\Ball_\Hil(1)} \Exs\bigl\| \widetilde{f}- f^*
\bigr\|_{L^2}^2. 
\end{equation}
Obtaining lower bounds on this minimax risk requires another
approximation property of the norm \mbox{$\|\cdot\|_{\Phi}$} relative to
\mbox{$\|\cdot\|_{L^{2}}$}. Consider matrix $\Psi\in\real^{m
\times m}$ with entries $\Psi_{ij}:=
\langle\psi_i, \psi_j \rangle_{\Phi}$. Since the eigenfunctions are
orthogonal in $L^2$, the deviation of $\Psi$ from the identity
measures how well the inner product defined by $\Phi$ approximates
the $L^2$-inner product over the first $m$ eigenfunctions of the
kernel operator. For proving lower bounds, we require an upper bound
of the form
\[
\mbox{(B2)} \quad\lambda_{\max}(\Psi) \leq c_{1}
\]
for some universal constant $c_{1}>0$. As the proof will clarify,
this upper bound is necessary in order that the Kullback--Leibler
divergence---which controls the relative discriminability between
different models---can be upper bounded in terms of the $L^2$-norm.
%
\begin{theorem}[(Lower bounds for $\|\widetilde{f}- f^*\|_{L^2}^2$)]
\label{THMMINIMAXLTWO} Suppose that condition \textup{(B2)} holds, and
the operator associated with kernel function $\KerFun$ of the
reproducing kernel Hilbert space $\Hil$ has eigenvalues with
polynomial-$\alpha$ decay.
\begin{longlist}[(a)]
\item[(a)] For the time-sampling model, the minimax risk
is lower bounded as
%
\begin{equation}
\label{eqminimaxL2time} \MiniMax^\Hil\bigl(\|\cdot\|_{L^2} \bigr)
\geq C \biggl\{ \min\biggl\{ \biggl(\frac{\sigma_0^2}{mn} \biggr
)^{{{2 \alpha}/({2 \alpha+1})}},
\frac{\sigma_0^2}{n} \biggr\} + \biggl(\frac{1}{m}\biggr)^{2 \alpha}
\biggr\}.
\end{equation}
\item[(b)] For the frequency-truncation model, the minimax error is
lower bounded as
%
\begin{equation}
\label{eqminimaxL2freq} \MiniMax^\Hil\bigl( \| \cdot\|_{L^2} \bigr)
\geq C \biggl\{ \biggl(\frac{\sigma_0^2}{ n} \biggr)^{{{2 \alpha
}/({2 \alpha+1})}} + \biggl(
\frac{1}{m} \biggr)^{2 \alpha} \biggr\}.
\end{equation}
\end{longlist}
\end{theorem}

Verifying condition (B2) requires, in general, some calculations
in the case of the time-sampling model. It is verified for uniform
time-sampling for the first-order Sobolev RKHS in
Appendix A 
of the supplementary material~\cite{supp}. For the frequency-truncation
model, condition (B2) always holds trivially since $\Psi=
I_{m}$. By this theorem, the $L^{2}$ convergence rates of
Corollaries~\ref{CorTimeFunc} and~\ref{CorBasisTrunFunc} are minimax
optimal. Also note that due to the presence of the approximation term
$m^{-2\alpha}$ in (\ref{eqminimaxL2freq}), the $\Phi$-norm
term $n^{{2\alpha}/({2\alpha+1})}$ is only dominant when
$m\ge c_{2} n^{{1}/({2\alpha+ 1})}$ implying that this
is the interesting regime for Theorem~\ref{THMMINIMAXEMP}(b).


\section{Proof of subspace-based rates}
\label{SecProofSubspace}

We now turn to the proofs of the results involving the error
$\dist_{\mathrm{HS}}(\widehat{\mathfrak{Z}},\mathfrak{Z}^*)$ between the
estimated $\widehat{\mathfrak{Z}}$ and
true subspace
$\mathfrak{Z}^*$. We begin by proving Theorem~\ref{THMSUBSPACEERR},
and then
turn to its corollaries.

\subsection{Preliminaries}

We begin with some preliminaries before proceeding to the heart of the
proof. Let us first introduce some convenient notation. Consider the
$n\times m$ matrices
\[
Y:=\lleft[\matrix{ \yobs{1} & \yobs{2} & \cdots&\yobs{n} }
\rright]^T \quad\mbox{and}\quad W:=\lleft[\matrix{ \wnoise{1} &
\wnoise{2} & \cdots& \wnoise{n} } \rright]^T,
\]
corresponding to the observation matrix $Y$ and noise matrix $W$,
respectively. In addition, we define the matrix $B:= (\beta_{ij})
\in\reals^{{n}\times r}$, and the diagonal matrix $S:=
\diag(s_1,\ldots,s_r) \in\reals^{r\times r}$. Recalling that
${Z^*}:= (z^*_1,\ldots,z^*_r) \in\reals^{m\times r}$, the
observation model (\ref{EqnLinObs}) can be written in the matrix form
$Y= B({Z^*}S)^T + \sigma_mW$. Moreover, let us define the
matrices $\widebar{B}:=\frac{B^T B}{{n}} \in\reals^{r\times r}$ and
$\widebar{W}:=\frac{W^T B}{{n}} \in\reals^{m\times r}$. Using this
notation, some algebra shows that the associated sample covariance
$\widehat{\Sigma}_{n}:=\frac{1}{{n}} Y^T Y$ can be written in the form
%
\begin{equation}
\label{EqnSamSplit} \widehat{\Sigma}_{n}= \underbrace{{Z^*}S\widebar
{B}S\bigl({Z^*}\bigr)^T}_{\Gamma} + \Delta_1 +
\Delta_2,
\end{equation}
where $\Delta_1:=\sigma_m [ \widebar{W}S({Z^*})^T + {Z^*}
S
\widebar{W}{}^T ]$ and $\Delta_2:=\sigma_m^2 \frac{W^T
W}{{n}}$.

Lemma~\ref{LemFeasible}, proved\vspace*{1pt} in Appendix B.3 
of the supplementary material~\cite{supp}, establishes the existence of a matrix
${\widetilde{Z}^*}\in{V_{r}(\reals^{m})}$ such that $\Range
({\widetilde{Z}^*})
= \Range({Z^*})$. As discussed earlier, due to the nature of the
Steifel manifold, there are many versions of this matrix ${\widetilde
{Z}^*}$, and
also of any optimal solution matrix $\widehat{Z}$, obtained via right
multiplication with an orthogonal matrix. For the subsequent
arguments, we need to work with a particular version of ${\widetilde
{Z}^*}$ (and
$\widehat{Z}$) that we describe here.

Let us fix some convenient versions of ${\widetilde{Z}^*}$ and
$\widehat{Z}$. As a
consequence of CS decomposition, as long as $r\leq m/2$,
there exist orthogonal matrices $U,V\in\reals^{r\times
r}$ and an orthogonal matrix $Q\in\reals^{m\times
m}$ such that
%
\begin{equation}
\label{eqCSdecompconseq} Q^T {\widetilde{Z}^*}U = %
\pmatrix{ I_r
\cr
0
\cr
0 } \quad\mbox{and}\quad Q^T
\widehat{Z}V = %
\pmatrix{ \widehat{C}
\cr
\widehat{S}
\cr
0},
\end{equation}
where $\widehat{C}= \diag(\widehat{c}_1,\ldots,\widehat{c}_r)$ and
$\widehat{S}=
\diag(\widehat{s}_1,\ldots,\widehat{s}_r)$ such that $1 \ge
\widehat{s}_1 \ge\cdots\ge
\widehat{s}_r\ge0$ and $\widehat{C}{}^2 + \widehat{S}{}^2 = I_r$. See
Bhatia~\cite{Bhatia1996}, Theorem VII.1.8, for details on this
decomposition. In the analysis to follow, we work with ${\widetilde
{Z}^*}U$ and
$\widehat{Z}V$ instead of ${\widetilde{Z}^*}$ and~$\widehat{Z}$. To
avoid extra
notation, from
now on, we will use ${\widetilde{Z}^*}$ and $\widehat{Z}$ for these new
versions, which
we refer to as \textit{properly aligned}. With this choice, we may
assume $U = V = I_r$ in the CS
decomposition (\ref{eqCSdecompconseq}).

The following lemma isolates some useful properties of properly
aligned subspaces:
%
\begin{lemma}
\label{LemManageable}
Let ${\widetilde{Z}^*}$ and $\widehat{Z}$ be properly aligned, and define
the matrices
%
\begin{equation}
\widehat{P}:={P}_{\widehat{Z}} - {P}_{{\widetilde{Z}^*}} = \widehat{Z}
\widehat{Z}{}^T - {\widetilde{Z}^*}\bigl({\widetilde{Z}^*}
\bigr)^T \quad\mbox{and}\quad\widehat{E}:=\widehat{Z}- {
\widetilde{Z}^*}.
\end{equation}
In terms of the CS decomposition (\ref{eqCSdecompconseq}), we have:
%
\begin{subequation}
%
\begin{eqnarray}
\label{EqnManageable}
\nnorm\widehat{E}\nnorm_{{\mathrm{HS}}} &\leq& \nnorm
\widehat{P}\nnorm_{{\mathrm{HS}}},
\\
\label{eqpertmidproof2}
\bigl({\widetilde{Z}^*}
\bigr)^T ( {P}_{\widetilde{Z}^*}- {P}_{\widehat{Z}}) {
\widetilde{Z}^*} &=& \widehat{S}{}^2 \quad\mbox{and}
\\
\label{eqpertmidproof3}
{d_{\mathrm{HS}}^2\bigl(
\widehat{Z},{\widetilde{Z}^*}\bigr)} &=& \nnorm{P}_{{\widetilde
{Z}^*}} -
{P}_{\widehat{Z}} \nnorm_{\mathrm{HS}}^2 \nonumber\\
&=& 2\nnnorml
\widehat{S}{}^2 \nnnormr_{\mathrm{HS}}^2 + 2\nnorm\widehat{C}
\widehat{S}\nnorm_{\mathrm{HS}}^2 \\
&=& 2 \sum
_{k} \widehat{s}_k^2\bigl(
\widehat{s}_k^2 + \widehat{c}_k^2
\bigr) = 2 \trace\bigl(\widehat{S}{}^2\bigr).\nonumber
\end{eqnarray}
\end{subequation}
\end{lemma}
\begin{pf}
From the CS decomposition (\ref{eqCSdecompconseq}), we have
\[
{\widetilde{Z}^*} \bigl({\widetilde{Z}^*}\bigr)^T - \widehat{Z}(
\widehat{Z})^T = Q \pmatrix{ \widehat{S}{}^2 & -\widehat{C}
\widehat{S}& 0
\cr
-\widehat{S}\widehat{C}& -\widehat{S}{}^2 & 0
\cr
0 &
0 & 0} Q^T,
\]
from which
relations (\ref{eqpertmidproof2}) and (\ref{eqpertmidproof3})
follow. From decomposition (\ref{eqCSdecompconseq}) and the proper
alignment condition $U = V = I_r$, we have
%
\begin{eqnarray}
\label{eqFrobprojdistineq} \nnorm\widehat{E}\nnorm_{\mathrm{HS}}^2
&=& \nnnorml Q^T\bigl(\widehat{Z}- {\widetilde{Z}^*}\bigr)
\nnnormr_{\mathrm{HS}}^2 = \nnorm I_r-\widehat{C}
\nnorm_{\mathrm{HS}}^2 + \nnorm\widehat{S}\nnorm_{\mathrm{HS}}^2
\nonumber\\[-8pt]\\[-8pt]
&=& 2 \sum_{i=1}^r(1 -
\widehat{c}_i) \le2 \sum_{i=1}^r
\bigl(1 - \widehat{c}_i^2\bigr) = 2 \sum
_{i=1}^r\widehat{s}_i^2 =
\nnorm\widehat{P}\nnorm_{\mathrm{HS}}^2,\nonumber
\end{eqnarray}
where we have used the relations $\widehat{C}{}^2 + \widehat{S}{}^2 =
I_r$, $\widehat{c}_i
\in[0,1]$ and $2 \trace(\widehat{S}{}^2) = \nnorm{P}_{{\widetilde
{Z}^*}} - {P}_{\widehat{Z}} \nnorm_{\mathrm{HS}}^2$.
\end{pf}


\subsection{\texorpdfstring{Proof of Theorem \protect\ref{THMSUBSPACEERR}}{Proof of Theorem 1}}

Using the notation introduced in Lemma~\ref{LemManageable}, our goal
is to bound the Hilbert--Schmidt norm $\nnorm\widehat{P}\nnorm_{{\mathrm{HS}}}$. Without
loss of
generality we will assume $s_1 = 1$ throughout. Recalling
definition (\ref{EqnSamSplit}) of the random matrix $\Delta$, the
following inequality plays a central role in the proof:

\begin{lemma}
\label{LemPerturbation}
Under condition \textup{(A1)} and $s_1 =1$, we have
%
\begin{equation}
\label{EqnBasic} \nnorm\widehat{P}\nnorm_{{\mathrm{HS}}}^2 \leq128
\llangle{\widehat{P}},{\Delta_1 + \Delta_2}\rrangle
\end{equation}
with probability at least $1 - \exp(-{n}/32)$.
\end{lemma}
\begin{pf}
We use the shorthand notation $\Delta= \Delta_1 + \Delta_2$ for the
proof. Since ${\widetilde{Z}^*}$ is feasible and $\widehat{Z}$ is optimal
for
problem (\ref{EqnDiscreteMest}), we have the basic inequality
$\llangle{\widehat{\Sigma}_{n}},{{P}_{\widetilde
{Z}^*}}\rrangle\leq
\llangle{\widehat{\Sigma}_{n}},{{P}_{\widehat{Z}}}\rrangle$. Using
the decomposition\vadjust{\goodbreak}
$\widehat{\Sigma}= \Gamma+ \Delta$ and rearranging yields the inequality
%
\begin{equation}
\label{eqpertmidproofinequ} \llangle{\Gamma},{{P}_{\widetilde{Z}^*}-
{P}_{\widehat{Z}}}\rrangle\leq\llangle{\Delta},{{P}_{\widehat{Z}}-
{P}_{\widetilde{Z}^*}}\rrangle.
\end{equation}
From definition (\ref{EqnSamSplit}) of $\Gamma$ and ${Z^*}=
{\widetilde{Z}^*}
R$, the left-hand side of the
inequality (\ref{eqpertmidproofinequ}) can be lower bounded as
\begin{eqnarray*}
\llangle{\Gamma},{{P}_{\widetilde{Z}^*}- {P}_{\widehat{Z}}}\rrangle
&=&
\llanglel{
\widebar{B}},{SR^T \bigl({\widetilde{Z}^*}\bigr)^T
({P}_{\widetilde{Z}^*}- {P}_{\widehat{Z}}) {\widetilde{Z}^*}R
S}\rrangler
\\
&=& \trace{\widebar{B}} {S R^T \widehat{S}{}^2
RS}
\\
&\geq& \lambda_{\min}(\widebar{B}) \lambda_{\min}
\bigl(S^2\bigr) \lambda_{\min
}\bigl(R^TR\bigr)
\trace\bigl( \widehat{S}{}^2\bigr),
\end{eqnarray*}
where we have used (90) 
and (91) 
of Appendix I several times 
(cf. the supplementary material~\cite{supp}). We note that
$\lambda_{\min}(S^2) = s_r^2 \ge\frac{1}{2}$ and
$\lambda_{\min}(R^TR) \ge\frac{1}{2}$ provided
$rC_m(f^*)\ge
\frac{1}{2}$; see equation (70). 
To bound the
minimum eigenvalue of $\widebar{B}$, let $\gamma_{\min}(B)$ denote the
minimum singular value of the ${n}\times r$ Gaussian matrix
$B$. The following concentration inequality is well known
(cf.~\cite{Davidson2001,Ledoux2001}):
\[
\mprob\bigl[ \gamma_{\min}(B) \leq\sqrt{{n}} - \sqrt{r} -t \bigr] \le
\exp\bigl(-t^2/2\bigr)\qquad\mbox{for all $t > 0$.}
\]
Since $\lambda_{\min}(\widebar{B}) = \gamma_{\min}^2(B/\sqrt{{n}})$, we
have that $\lambda_{\min}(\widebar{B}) \ge(1 - \sqrt{r/{n}} - t)^2$
with probability at least $1-\exp(-{n}t^2 /2)$. Assuming $r
/{n}\le\frac{1}{4}$ and setting $t = \frac{1}{4}$, we get
$\lambda_{\min}(\widebar{B}) \ge\frac{1}{16}$ with probability at least
$1-\exp(-{n}/32)$. Putting the pieces together yields the claim.
\end{pf}

Inequality (\ref{EqnBasic}) reduces the problem of bounding
$\nnorm\widehat{P}\nnorm_{{\mathrm{HS}}}^2$ to the sub-problem of studying
the random
variable $\llangle{\widehat{P}},{\Delta_1 + \Delta_2}\rrangle$.
Based on
Lemma~\ref{LemPerturbation}, our next step is to establish an
inequality (holding with high probability) of the form
%
\begin{equation}
\label{EqnInitialBound} \llangle{\widehat{P}},{\Delta_1 +
\Delta_2}\rrangle\leq c_1 \biggl\{ \frac{\sigma_m}{\sqrt{{n}}}
r^{3/2} \SimpleF\bigl( \nnorm\widehat{E}\nnorm_{\mathrm{HS}}\bigr) + \kappa\nnorm
\widehat{E}\nnorm_{\mathrm{HS}}^2 + \epscrit^2 \biggr\},
\end{equation}
where $c_1$ is some universal constant, $\kappa$ is the constant in
condition (A1) and $\epscrit$ is the critical radius from
Theorem~\ref{THMSUBSPACEERR}. Doing so is a nontrivial task: both
matrices $\widehat{P}$ and $\Delta$ are random and depend on one another,
since the subspace $\widehat{Z}$ was obtained by optimizing a random function
depending on $\Delta$. Consequently, our proof of
bound (\ref{EqnInitialBound}) involves deriving a uniform law of large
numbers for a certain matrix class.

Suppose that bound (\ref{EqnInitialBound}) holds, and that the
subspaces ${\widetilde{Z}^*}$ and $\widehat{Z}$ are properly aligned.
Lemma~\ref{LemManageable} implies that $\nnorm\widehat{E}\nnorm_{{\mathrm{HS}}}
\leq
\nnorm\widehat{P}\nnorm_{{\mathrm{HS}}}$, and since $\SimpleF$ is a nondecreasing
function, inequality (\ref{EqnInitialBound}) combined with
Lemma~\ref{LemPerturbation} implies that
\[
(1 - 128 \kappa c_1 ) \nnorm\widehat{P}\nnorm_{{\mathrm{HS}}}^2
\leq c_1 \biggl\{ \frac{\sigma_m}{\sqrt{{n}}} r^{3/2} \SimpleF\bigl(
\nnorm\widehat{P}\nnorm_{\mathrm{HS}}\bigr) + \epscrit^2 \biggr\},
\]
from which\vspace*{2pt} the claim follows as long as $\kappa$ is suitably small
(e.g., $\kappa\leq\frac{1}{256 c_1}$ suffices). Accordingly,
in order to complete the proof of Theorem~\ref{THMSUBSPACEERR}, it\vadjust{\goodbreak}
remains to prove bound (\ref{EqnInitialBound}), and the remainder
of our work is devoted to this goal. Given the linearity of trace, we
can bound the terms $\llangle{\widehat{P}},{\Delta_1}\rrangle$ and
$\llangle{\widehat{P}},{\Delta_2}\rrangle$ separately.

\subsubsection{\texorpdfstring{Bounding $\llangle{\widehat{P}},{\Delta_1}\rrangle$}
{Bounding <<P, Delta1>>}}

Let $\{\overline{z}_j\}$, $\{\widetilde{z}^*_j\}$ and $\{\widehat
{e}_j\}$ and $\{\widebar{w}_j\}$ denote
the columns of $\widehat{Z}$, ${\widetilde{Z}^*}$, $\widehat{E}$ and
$\widebar{W}$,
respectively, where
we recall the definitions of these quantities from
equation (\ref{EqnSamSplit}) and Lemma~\ref{LemManageable}. Note that
$\widebar{w}_j = {n}^{-1} \sum_{i=1}^{n}\wnoise{i} \beta_{ij}$.
In Appendix C.1 
of the supplementary material~\cite{supp}, we show that
%
\begin{equation}
\label{EqnDeltaOne} \llangle{\widehat{P}},{\Delta_1}\rrangle\leq
\sqrt{6} \sigma r^{3/2} \max_{j,k} \bigl|\langle
\widebar{w}_k,\widehat{e}_j\rangle\bigr| + \sqrt{
\frac{3}{2}} \sigma r \nnorm\widehat{E}\nnorm_{\mathrm{HS}}^2
\max_{j,k} \bigl|\bigl\langle\widebar{w}_j,
\widetilde{z}^*_k\bigr\rangle\bigr|.
\end{equation}
Consequently, we need to obtain bounds on quantities of the form
$|\langle\widebar{w}_j,v\rangle|$, where the vector $v$ is either fixed
(e.g., $v
= \widetilde{z}^*_j$) or random (e.g., $v = \widehat{e}_j$). The
following lemmas
provide us with the requisite bounds:
%
\begin{lemma}
\label{LemDeltaOneEllipse}
We have
\[
\max_{j,k} \sigma r^{3/2} \bigl| \langle\widebar{w}_k,
\widehat{e}_j\rangle\bigr| \le C \biggl\{ \frac{\sigma}{\sqrt{{n}}}
r^{3/2} \SimpleF\bigl(\nnorm\widehat{E}\nnorm_{\mathrm{HS}}\bigr) + \kappa\nnorm
\widehat{E}\nnorm_{\mathrm{HS}}^2 + \kappa\epscrit^2
\biggr\}
\]
with probability at least $1-c_1 r\exp(-\kappa^2 r^{-3}
{n}
\frac{\epscrit^2}{2\sigma^2}) - r\exp(-{n}/64)$.
\end{lemma}
%
%
\begin{lemma}
\label{LemDeltaOneFixed}
We have
\[
\mprob\Bigl[\max_{j,k} \sigma r\bigl| \widebar{w}_k^T
\widetilde{z}^*_j\bigr| \leq\sqrt{6} \kappa\Bigr] \geq1 -
r^2 \exp\bigl(- \kappa^2 r^{-2} {n}/2
\sigma^2\bigr).
\]
\end{lemma}

See Appendices C.2 
and C.3 
in the supplementary material~\cite{supp} for the proofs of these
claims.


\subsubsection{\texorpdfstring{Bounding $\llangle{\widehat{P}},{\Delta_2}\rrangle$}
{Bounding <<P, Delta2>>}}

Recalling definition (\ref{EqnSamSplit}) of $\Delta_2$ and using
linearity of the trace, we obtain
\[
\llangle{\widehat{P}},{\Delta_2}\rrangle= \frac{\sigma^2}{{n}} \sum
_{j=1}^r \bigl\{ (\overline{z}_j)^T
W^T W\overline{z}_j - \bigl(\widetilde
{z}^*_j\bigr)^T W^T W\widetilde{z}^*_j
\bigr\}.
\]
Since $\widehat{e}_j = \overline{z}_j - \widetilde{z}^*_j$, we have
%
\begin{eqnarray}\label{eqGDelta2expansion}
\llangle{\widehat{P}},{\Delta_2}\rrangle& = & \sigma^2
\sum_{j=1}^r \biggl\{ 2 \bigl(
\widetilde{z}^*_j\bigr)^T \biggl( \frac{1}{{n}}W^TW-
I_r \biggr) \widehat{e}_j + \frac{1}{{n}} \| W
\widehat{e}_j\|_2^2 + 2 \bigl(
\widetilde{z}^*_j\bigr)^T \widehat{e}_j
\biggr\}
\nonumber\\[-8pt]\\[-8pt]
& \leq &
\sigma^2 \sum_{j=1}^r \biggl\{
2 \underbrace{\bigl(\widetilde{z}^*_j\bigr)^T \biggl(
\frac{1}{{n}}W^TW- I_r \biggr) \widehat{e}_j}_{T_1(\widehat{e}_j;
\widetilde{z}^*_j)}
+ \underbrace{\frac{1}{{n}} \| W\widehat{e}_j
\|_2^2}_{T_2(\widehat{e}_j)} \biggr\},\nonumber
\end{eqnarray}
where we have used the fact that $2\sum_j (\widetilde{z}^*_j)^T
\widehat{e}_j = 2\sum_j [
(\widetilde{z}^*_j)^T\overline{z}_j - 1] = 2\sum_j (\widehat{c}_j -
1) = -\nnorm
\widehat{E}\nnorm_{\mathrm{HS}}^2 \le
0$.

The following lemmas provide high probability bounds on the terms
$T_1$ and $T_2$.
%
\begin{lemma}
\label{LemTermOne}
We have the upper bound
\[
\sigma^2 \sum_{j=1}^rT_1
\bigl(\widehat{e}_j;\widetilde{z}^*_j\bigr) \le C \biggl
\{ \sigma_0 \frac{\sigma}{\sqrt{{n}}} r \SimpleF\bigl(\nnorm\widehat{E}
\nnorm_{\mathrm{HS}}\bigr) + \kappa\nnorm\widehat{E}\nnorm_{\mathrm{HS}}^2
+ \kappa\epscrit^2 \biggr\}
\]
with probability $1-c_2 \exp(- \kappa^2 r^{-2} {n} \frac
{\epscrit\wedge\epscrit^2}{16 \sigma^2}) - r\exp(-{n}/64)$.
\end{lemma}
%
\begin{lemma}
\label{LemTermTwo}
We have the upper bound $\sigma^2 \sum_{j=1}^rT_2(\widehat{e}_j) \le
C
\kappa\{ \nnorm\widehat{E}\nnorm_{\mathrm{HS}}^2 + \epscrit^2 \}$ with
probability at least $1-c_3 \exp(-\kappa^2 r^{-2} {n}
\epscrit^2/2\sigma^2)$.
\end{lemma}

See Appendices C.4 
and C.5 
in the supplementary material~\cite{supp} for the proofs of these claims.


\section{Discussion}\label{SecDiscuss}

We studied the problem of sampling for functional PCA from a
functional-theoretic viewpoint. The principal components were assumed
to lie in some Hilbert subspace $\Hil$ of $L^{2}$, usually a RKHS, and
the sampling operator, a bounded linear map $\Phi\dvtx  \Hil\to
\reals^{m}$. The observation model was taken to be the output of
$\Phi$ plus some Gaussian noise. The two main examples of $\Phi$
considered were time sampling, $[\Phi f]_{j} = f(t_{j})$ and
(generalized) frequency truncation $[\Phi f]_{j} =
\langle\psi_{j}, f \rangle_{L^{2}}$. We showed that it is possible
to recover
the subspace spanned by the original components, by applying a regularized
version of PCA in $\real^{m}$ followed by simple linear mapping
back to function space. The regularization involved the
``trace-smoothness condition'' (\ref{EqnTraceSmooth}) based on the
matrix $K= \Phi\Phi^{*}$ whose eigendecay influenced the rate
of convergence in $\reals^{m}$.

We obtained the rates of convergence for the subspace estimators both
in the discrete domain, $\reals^{m}$, and the function domain,
$L^{2}$. As examples, for the case of a RKHS $\Hil$ for which both
the kernel integral operator and the kernel matrix $K$ have
polynomial-$\alpha$ eigendecay (i.e., $\mu_{j} \asymp\widehat{\mu}_{j}
\asymp j^{-2\alpha}$), the following rates in $\mathrm{HS}$-projection
distance for subspaces in the function domain were worked out in
detail:\vspace*{10pt}

\begin{center}
\begin{tabular}[h]{@{}cc@{}}
\hline
Time sampling & Frequency truncation \\
\hline
$ ( \frac{1}{mn} )^{{2 \alpha}/({2\alpha+
1})} + ( \frac{1}{m} )^{2\alpha}$ &
$ ( \frac{1}{n} )^{{2 \alpha}/({2\alpha+
1})} + ( \frac{1}{m} )^{2\alpha}$\\
\hline\\
\end{tabular}
\end{center}

\noindent The two terms in each rate can be associated, respectively, with the
estimation error (due to noise) and approximation error (due to having
finite samples of an infinite-dimensional object). Both rates exhibit a
trade-off between the number of statistical samples\vadjust{\goodbreak} ($n$) and that of
functional samples ($m$). The two rates are qualitatively different:
the two terms in the time sampling case interact to give an overall
fast rate of $n^{-1}$ for the optimal trade-off $m \asymp
n^{1/({2\alpha})}$, while there is no interaction between the two
terms in the frequency truncation; the optimal trade-off gives an
overall rate of $n^{-{2\alpha }/({2\alpha+1 })}$, a characteristics
of nonparametric problems. Finally, these rates were shown to be
minimax optimal.

\begin{supplement}
\stitle{Proofs and auxiliary results}
\slink[doi]{10.1214/12-AOS1033SUPP} 
\sdatatype{.pdf}
\sfilename{aos1033\_supp.pdf}
\sdescription{This supplement contains some of the proofs and
auxiliary results referenced in the text.}
\end{supplement}


\printaddresses

\end{document}